\documentclass[10pt,reqno,a4paper]{amsart}

\usepackage[T1]{fontenc}
\usepackage[utf8]{inputenc}
\usepackage{lmodern}
\usepackage{amsmath,amssymb,mathtools,mathrsfs}
\usepackage{bm}
\usepackage{enumitem}
\usepackage{microtype}
\usepackage[a4paper,top=24mm,bottom=25mm,left=25mm,right=25mm,headheight=13pt,headsep=7mm,footskip=10mm]{geometry}
\usepackage{graphicx}
\usepackage{float}
\usepackage{needspace}
\usepackage{placeins}
\usepackage[colorlinks=true,linkcolor=blue,citecolor=blue,urlcolor=blue,hypertexnames=false]{hyperref}

\allowdisplaybreaks
\numberwithin{equation}{section}

\newtheorem{theorem}{Theorem}[section]
\newtheorem{proposition}[theorem]{Proposition}
\newtheorem{lemma}[theorem]{Lemma}
\newtheorem{corollary}[theorem]{Corollary}
\theoremstyle{definition}
\newtheorem{assumption}[theorem]{Assumption}
\newtheorem{definition}[theorem]{Definition}
\newtheorem{example}[theorem]{Example}
\theoremstyle{remark}
\newtheorem{remark}[theorem]{Remark}

\newcommand{\R}{\mathbb R}
\newcommand{\C}{\mathbb C}
\newcommand{\Real}{\operatorname{Re}}

\title[Two-sided weighted-history transfer]{Two-Sided Weighted-History Transfer for Abstract Systems with Infinite Memory}

\author[G.~J.-M.~Bomisso]{Gossrin Jean-Marc Bomisso}
\makeatletter
\g@addto@macro\@setauthors{\par{\centering\normalfont\small Laboratoire de Math\'ematiques et Informatique, UFR Sciences Fondamentales et Appliqu\'ees\\ Universit\'e Nangui Abrogoua, Abidjan, C\^ote d'Ivoire\\ \texttt{bogojm@yahoo.fr}\par}}
\makeatother

\makeatletter
\@ifundefined{subjclassname@2020}{\@namedef{subjclassname@2020}{2020 Mathematics Subject Classification}}{}
\makeatother

\begin{document}

\begin{abstract}
We prove a two-sided transfer principle for second-order systems with Hilbert-valued infinite memory. The squared norm of the coupled semigroup, from a centered weighted-history space into the energy space, is bounded above and below by constant multiples of $\rho_w(t)$. Here $\rho_w(t)$ is the squared norm of the right-translation semigroup between the corresponding weighted and unweighted memory spaces. The coupled system satisfies an energy identity, a mechanical multiplier estimate and a bounded observation. The memory density has decay of order at most polynomial, and the weight is admissible. Under these hypotheses the energy is bounded by the initial norm times $\kappa(t)+\rho_w(t)$, where $\kappa$ is the tail of the kernel. The two-sided equivalence holds when $\rho_w$ is integrable and the kernel has a finite first moment and satisfies a tail bound, or when $\kappa\le C_\kappa\rho_w$. For the density $pb(1+s)^{-p-1}$, with $p>1$ and $b>0$, and the weight $(1+s)^{2\theta}$, with $\theta>0$, the rate is $(1+t)^{-\min\{2\theta,p+1\}}$. For kernels and weights with logarithmic factors, under explicit conditions on the exponents, the logarithm of the weight enters the rate below the saturation threshold, and the logarithm of the kernel enters above it. A theorem with a prescribed decay of the inherited history gives exact exponents for algebraically growing histories. The hypotheses are verified for gradient-memory waves without instantaneous damping, endpoint-damped strings, damped plates and oscillators.
\end{abstract}

\maketitle

\section{Introduction}
\label{sec:introduction}

A material with memory responds not only to its present deformation but also to its past motion. In the history formulation of Dafermos \cite{Dafermos1970}, this dependence becomes a transport equation coupled to the mechanical evolution. The corresponding semigroup acts on an augmented energy space, within the standard theory of strongly continuous semigroups \cite{Pazy}. Its history component permits initial data extending arbitrarily far into the past. Such data can retain energy for a long time, even when the memory kernel is integrable and the mechanical dynamics are dissipative. The formulation therefore raises a quantitative question: \emph{how does the localization of the initial history determine the rate at which the coupled system loses energy?}

Two quantitative questions arise for prescribed histories. One asks for the decay of a single solution from its given past. The other seeks an estimate uniform over the unit ball of an initial-data space. The latter formulation makes it possible to compare decay rates as the localization of the initial history is strengthened. Weighted spaces provide a natural scale for this comparison. We determine the corresponding operator-norm rates and their optimal exponents.

Decay estimates for equations with memory have a substantial abstract literature. Alabau-Boussouira, Cannarsa and Sforza \cite{AlabauCannarsaSforza2008} developed integral inequalities and multiplier techniques for second-order evolution equations. Fabrizio and Polidoro \cite{FabrizioPolidoro2002} showed that exponential decay requires an exponentially decaying kernel. Polynomial rates for weakly dissipative second-order memory systems were obtained by Mu\~noz Rivera, Naso and Vegni \cite{MunozRiveraNasoVegni2003} and by Mu\~noz Rivera and Naso \cite{MunozRiveraNaso2006,MunozRiveraNaso2010}. General decay estimates under broad assumptions on the kernel were developed by Guesmia \cite{Guesmia2011,Guesmia2015,Guesmia2020}, Guesmia and Messaoudi \cite{GuesmiaMessaoudi2014}, Mustafa \cite{Mustafa2021}, and Jin, Liang and Xiao \cite{JinLiangXiao2014,JinLiangXiao2019,JinLiangXiao2023}. These works give upper bounds for trajectories or for classes of data.

The data spaces of these results differ from the one considered here. Guesmia and Messaoudi \cite[Theorem~2.1]{GuesmiaMessaoudi2014} assume a bounded past history, and their bound contains the tail $\int_t^\infty g$ of the kernel. Guesmia \cite[Theorem~3.1]{Guesmia2015} treats coupled systems with indirect damping, and the rate is governed by the regularity of the data. Jin, Liang and Xiao \cite[Theorem~3.1]{JinLiangXiao2019} prove the rate $t^{-1}$ for every nonincreasing integrable kernel, for data of higher regularity whose histories are weighted by the tail of the kernel. Conti and Pata \cite[Theorem~4.4 and Corollary~4.7]{ContiPata2020} obtain decay for the Volterra equation without prehistory and kernels with flat zones and jumps. Their Proposition~1.1 shows that the energy of each nonzero solution cannot decay faster than the kernel. Proposition~4.4 of Conti, Dell'Oro and Pata \cite{ContiDellOroPata2022} bounds the full semigroup from below by the history translation on the natural memory space, and their Theorem~7.3 gives a lower rate for individual trajectories. Dell'Oro, Li and Pata \cite{DellOroLiPata2026} study the translation generator and its spectrum. For a boundary-memory Euler--Bernoulli beam, Anass\'e, Bomisso and Kpata \cite{AnasseBomissoKpata2026} established optimal polynomial decay for compatible trajectories with zero prehistory.

The present result is an equivalence of operator norms. The initial history is measured, after centering, in a weighted space $\mathcal M_w$. We prove that the squared norm of the coupled semigroup from this weighted space into the energy space is comparable, from above and from below, to the exact squared norm $\rho_w(t)$ of the weighted right translation. The rate is thus computed by an explicit supremum built from the kernel and the weight. The lower bound extends the translation mechanism of \cite{ContiDellOroPata2022} to weighted spaces. The upper bound, uniform over the weighted unit ball and without regularity of the data, is the main new estimate. This framework places several earlier rates in a common setting. Bounded histories produce the tail term of \cite{GuesmiaMessaoudi2014}. The weight $1+\kappa/m$ built from the tail $\kappa$ of the kernel yields the universal upper bound $\rho_w(t)\le2/(1+t)$ for every nonincreasing integrable kernel. Under the structural assumptions of Corollary~\ref{cor:tail-weight}, it gives an exact two-sided semigroup law of order $(1+t)^{-1}$, which recovers the rate of \cite{JinLiangXiao2019}.

We work on a Gelfand triple $V\hookrightarrow H\hookrightarrow V'$ and allow the memory observation to take values in a Hilbert space $Y$. The mechanical part is described by a coercive elastic form and a dissipative form. A bounded map $\ell:V\to Y$ couples it to a relative-history variable $z$ with memory density $m$. We measure the initial history after subtracting the initial observed displacement. Thus $\xi_0=z_0-\ell(u_0)$ is the centered history. This centering follows from the characteristic representation of the transport equation. It separates the inherited past from the history generated after the initial time. The energy space is denoted by $\mathscr X$. For a weight $w\ge1$, the space $\mathscr X_w$ consists of energy states whose centered histories belong to $\mathcal M_w=L^2(\R_+,m(s)w(s)\,ds;Y)$. The squared norm of the right translation from $\mathcal M_w$ into $\mathcal M_m=L^2(\R_+,m(s)\,ds;Y)$ is given by
\[
 \rho_w(t)=\operatorname*{ess\,sup}_{r\ge0}\frac{m(r+t)}{m(r)w(r)} .
\]

Our main result is the transfer principle of Theorem~\ref{thm:general-transfer} and Corollary~\ref{cor:general-two-sided}. Assume that $m$ has decay of order at most polynomial and satisfies a logarithmic-derivative bound, and that $w$ is admissible. Then the energy is bounded by the initial norm in $\mathscr X_w$ times $\kappa(t)+\rho_w(t)$, where $\kappa$ is the tail of $m$. If $\rho_w$ is integrable and $m$ has a finite first moment and satisfies a tail bound, or if $\kappa\le C_\kappa\rho_w$, there are constants $0<c\le C$ such that
\begin{equation*}
 c\,\rho_w(t)\le\|T(t)\|^2_{\mathcal L(\mathscr X_w,\mathscr X)}\le C\rho_w(t),\quad\text{for all }t\ge0 .
\end{equation*}
Both bounds are uniform in the initial state. The rate of the coupled system is therefore the rate of the weighted right translation. Individual trajectories may decay faster.

All results of the paper are applications of one theorem. Theorem~\ref{thm:source-transfer} converts any envelope $K\sigma(t)$ of the inherited energy into the energy bound $(E(0)+K)(\kappa(t)+\sigma(t))$, and into $(E(0)+K)(m(t)+\sigma(t))$ when $\sigma$ is integrable. The weighted transfer is the case $\sigma=\rho_w$. The prescribed-history theorem is the case $\sigma(t)=(1+t)^{-\beta}$. The polynomial and logarithmic laws follow by computing $\rho_w$.

The polynomial case is a fundamental consequence. Let the system satisfy Assumption~\ref{ass:decay-transfer}. For $m(s)=pb(1+s)^{-p-1}$, with $p>1$ and $b>0$, and $w(s)=(1+s)^{2\theta}$, with $\theta>0$, Theorem~\ref{thm:sharp-abstract-transfer} gives positive constants $c_{p,\theta}$ and $C_{p,\theta}$ such that
\begin{equation*}
 c_{p,\theta}(1+t)^{-\alpha(p,\theta)}
 \le \|T(t)\|_{\mathcal L(\mathscr X_\theta,\mathscr X)}^2
 \le C_{p,\theta}(1+t)^{-\alpha(p,\theta)}
 \quad\text{with}\quad \alpha(p,\theta)=\min\{2\theta,p+1\}.
\end{equation*}
The rate improves with the localization parameter until $\theta=(p+1)/2$. Beyond this threshold, the exponent remains $p+1$. This case is deduced from Corollary~\ref{cor:general-two-sided}.

Two further consequences show the role of the pair $(m,w)$. For the tail weight $w=1+\kappa/m$, Corollary~\ref{cor:tail-weight} gives the exact rate $(1+t)^{-1}$, under (K3) and a finite first moment of $m$. For a kernel $(1+s)^{-p-1}\log(e+s)^{\alpha}$ and a weight $(1+s)^{2\theta}\log(e+s)^{\eta}$, Theorem~\ref{thm:log-rates} computes $\rho_w$ with explicit constants, for $p>1$, $\eta\ge0$ and $\alpha$ satisfying \eqref{eq:log-alpha}. Below the saturation threshold the logarithm of the weight enters the rate, and above it the logarithm of the kernel. Corollary~\ref{cor:log-waves} gives the corresponding two-sided law for viscoelastic waves without instantaneous damping.

For polynomial kernels, Theorem~\ref{thm:prescribed-source-transfer} addresses histories described by their inherited energy rather than a weighted norm. It is the case $\sigma(t)=(1+t)^{-\beta}$ of Theorem~\ref{thm:source-transfer}. If the translated history tail is bounded by $K(1+t)^{-\beta}$, with $\beta>0$, the coupled energy decays with exponent $\min\{\beta,p+1\}$. Histories with matching two-sided algebraic growth lead to sharp exponents for individual trajectories. Applied to the polynomial-kernel subclass of the first wave model in \cite{Guesmia2020}, this argument removes the logarithmic factor in the borderline upper estimate of equation~(4.28) of \cite{Guesmia2020}, at $q-r=2$. For nonzero compactly supported pure histories, the two-sided decay estimate also holds along each fixed trajectory, as shown in Section~\ref{sec:other-applications}.

We distinguish the result from the method of proof. The result is the transfer principle, namely the two-sided comparison of the two operator norms. The method is what we call the \emph{weighted right-translation norm approach}. It keeps the inherited history as a source term bounded by an envelope $\sigma$, equal to $\rho_w$ for weighted data, throughout a two-pass memory--dissipation splitting. The first pass gives an integrable observation, and the second reaches the rate of the kernel. A two-pass argument of this type was used in \cite{AnasseBomissoKpata2026} for zero-prehistory beam trajectories, where the source vanishes. The splitting age is chosen proportional to the dyadic time, which leads to a linear dyadic recursion (Lemmas~\ref{lem:linear-recursion} and \ref{lem:dyadic-step}). The lower bound follows directly from the characteristic formula. We verify the required mechanical estimate both when instantaneous damping controls the motion and when a coercive memory observation supplies the missing control. The latter case includes gradient-memory waves without instantaneous damping. Other examples concern a string damped at one endpoint, a damped plate and an oscillator.

The paper is organized as follows. Section~\ref{sec:variational} formulates the abstract system and its generator, develops the multiplier estimates and proves the transfer principle. It then deduces the explicit exponents for polynomial kernels and proves the prescribed-source theorem. Section~\ref{sec:other-applications} treats oscillators, strings, plates, gradient-memory waves, and logarithmic kernels and weights. It also presents a numerical illustration. Section~\ref{sec:conclusion} concludes.

\section{Abstract framework and the weighted-history transfer}
\label{sec:variational}

We consider a Hilbert-valued memory observation. Scalar boundary memory corresponds to $Y=\C$.

\subsection{Gelfand triple, mechanical forms and history observation}
\label{sec:variational-setting}

Let $V$ and $H$ be complex Hilbert spaces and assume a continuous, dense
embedding $V\hookrightarrow H$. We identify $H$ with its anti-dual, the space of continuous conjugate-linear functionals, so that
$V\hookrightarrow H\hookrightarrow V'$ is a Gelfand triple. All Hilbert
inner products are linear in the first argument. The symbol
$\langle\cdot,\cdot\rangle_{V',V}$ denotes the duality between the
anti-dual $V'$ and $V$, also linear in its first argument. We write $\R_+=[0,\infty)$ and denote by $\mathcal L(E,F)$ the space of bounded linear operators from $E$ to $F$. For a set $A$, $\mathbf 1_A$ denotes its characteristic function. Throughout, $C$ denotes a positive constant that may change from one estimate to another. Let
\begin{equation*}
 a:V\times V\to\C\quad\text{and}\quad d:V\times V\to\C
\end{equation*}
be continuous sesquilinear forms. We assume that $a$ is Hermitian and
$V$-coercive, and that $d$ is accretive:
\begin{equation}
 a(u,u)\ge c_a\|u\|_V^2\quad\text{and}\quad
 \Real d(v,v)\ge0,\quad\text{for all }u,v\in V.
 \label{eq:variational-basic-assumptions}
\end{equation}
The damping form need not be Hermitian. Let $Y$ be a nonzero complex Hilbert space and let
$\ell\in\mathcal L(V,Y)$ be a nonzero bounded linear operator. It may
represent a trace, a finite-dimensional observation, or a spatial
gradient. Let $\gamma>0$.

For the generation result, take a positive, locally absolutely continuous,
nonincreasing weight $m\in L^1(\R_+)$ and assume $m'$ exists almost
everywhere. The memory space is the weighted space $\mathcal M_m=L_m^2(\R_+;Y)$, endowed with the norm defined by
\[
 \|z\|_{\mathcal M_m}^2=\int_0^\infty m(s)\|z(s)\|_Y^2\,ds .
\]
The memory moment and the total mass of the kernel are given by
\[
 Jz=\int_0^\infty m(s)z(s)\,ds\quad\text{and}\quad\kappa_0=\int_0^\infty m(s)\,ds .
\]
The Cauchy--Schwarz inequality gives $\|Jz\|_Y\le\kappa_0^{1/2}\|z\|_{\mathcal M_m}$, so that $Jz$ is a Bochner integral in $Y$.
The mechanical system is the abstract variational equation given by
\begin{equation}
 \left\{\begin{aligned}
 \Bigl\langle \frac{\partial^2u}{\partial t^2}(t),\psi\Bigr\rangle_{V',V}+a(u(t),\psi)
  +d\Bigl(\frac{\partial u}{\partial t}(t),\psi\Bigr)+\gamma (Jz(t),\ell(\psi))_Y&=0,\quad\text{for all }\psi\in V,\\
 \frac{\partial z}{\partial t}(t,s)+\frac{\partial z}{\partial s}(t,s)&=\ell\Bigl(\frac{\partial u}{\partial t}(t)\Bigr),\quad\text{for all }s>0,\\
 z(t,0)&=0.
 \end{aligned}\right.
 \label{eq:variational-system}
\end{equation}
The $Y$-valued source $\ell(u_t)$ is meaningful for classical solutions with
$u_t\in V$. Mild solutions are defined through the generator below.

Set $\mathscr X=V\times H\times\mathcal M_m$ with the energy inner product defined by
\begin{equation*}
 \langle (u,v,z),(\widetilde u,\widetilde v,\widetilde z)\rangle_{\mathscr X}
 =a(u,\widetilde u)+(v,\widetilde v)_H
 +\gamma\int_0^\infty m(s)(z(s),\widetilde z(s))_Y\,ds.
\end{equation*}
Denote by $A_a,\mathcal D_d:V\to V'$ the form operators defined by
$\langle A_au,\psi\rangle=a(u,\psi)$ and
$\langle\mathcal D_dv,\psi\rangle=d(v,\psi)$, and by $\ell^*: Y\to V'$
the map $\langle\ell^*y,\psi\rangle=(y,\ell(\psi))_Y$.
The associated operator is defined as follows:
\begin{equation*}
 \mathscr A_m(u,v,z)
 =\bigl(v,-A_au-\mathcal D_dv-\gamma\ell^*Jz,\ell(v)-z_s\bigr)
\end{equation*}
with the domain given by
\begin{equation}
 \begin{split}
 D(\mathscr A_m)=\{(u,v,z)\in\mathscr X:\;&v\in V,\quad
 z\in D(\mathcal G),\\
 &A_au+\mathcal D_dv+\gamma\ell^*Jz\in H\subset V'\}.
 \end{split}
 \label{eq:variational-domain}
\end{equation}
Here $\mathcal G=-\partial_s$ is the generator of the
weighted right-translation semigroup on $\mathcal M_m$, with domain given by
\[
 D(\mathcal G)=\{z\in\mathcal M_m:
 z\in W^{1,1}_{\mathrm{loc}}((0,\infty);Y),\ 
 z_s\in\mathcal M_m,\ z(0)=0\}.
\]
The derivative $z_s$ is taken in the weak sense and $z(0)$ is its
$Y$-valued trace. The combined condition in the second line of
\eqref{eq:variational-domain} incorporates the mechanical operator,
damping and memory feedback. It also accounts for boundary forces
in the variational formulation.

\begin{lemma}[Weighted integration by parts]
\label{lem:weighted-ibp}
Let $E$ be a Hilbert space, and let $m$ be positive, locally absolutely continuous on $\R_+$, nonincreasing and integrable. Let $\eta\in W^{1,1}_{\rm loc}((0,\infty);E)$ satisfy $\eta,\eta_s\in L^2_m(\R_+;E)$ and $\eta(0)=0$. Then the following properties hold.
\begin{enumerate}[label=\textup{(\roman*)},leftmargin=*]
\item The following estimate and identity hold:
\[
 \int_0^\infty[-m'(s)]\|\eta(s)\|_E^2\,ds\le2\|\eta_s\|_{L^2_m(E)}\|\eta\|_{L^2_m(E)},
\]
\[
 \Real\int_0^\infty m(s)(\eta_s(s),\eta(s))_E\,ds=\frac12\int_0^\infty[-m'(s)]\|\eta(s)\|_E^2\,ds .
\]
\item $m(R)\|\eta(R)\|_E^2\to0$ as $R\to\infty$.
\item The integral $\int_0^\infty[-m'(s)]\eta(s)\,ds$ converges absolutely in $E$, with
\[
 \int_0^\infty[-m'(s)]\|\eta(s)\|_E\,ds\le m(0)^{1/2}\Bigl(\int_0^\infty[-m'(s)]\|\eta(s)\|_E^2\,ds\Bigr)^{1/2},
\]
and
\[
 \int_0^\infty m(s)\eta_s(s)\,ds=\int_0^\infty[-m'(s)]\eta(s)\,ds .
\]
\end{enumerate}
\end{lemma}
\begin{proof}
Fix $R>0$. Since $m\ge m(R)>0$ on $[0,R]$, the Cauchy--Schwarz inequality gives $\eta_s\in L^1(0,R;E)$. Hence $\eta$ is absolutely continuous on $[0,R]$. The functions $\|\eta\|_E^2$, $m\|\eta\|_E^2$ and $m\eta$ are then absolutely continuous on $[0,R]$, with
\[
 \bigl(m\|\eta\|_E^2\bigr)'=m'\|\eta\|_E^2+2m\Real(\eta_s,\eta)_E\quad\text{and}\quad(m\eta)'=m'\eta+m\eta_s .
\]
Integration on $[0,R]$ and $\eta(0)=0$ give
\begin{equation}
 m(R)\|\eta(R)\|_E^2+\int_0^R[-m'(s)]\|\eta(s)\|_E^2\,ds=2\Real\int_0^Rm(s)(\eta_s(s),\eta(s))_E\,ds .
 \label{eq:weighted-ibp-R}
\end{equation}
The right-hand side is bounded by $2\|\eta_s\|_{L^2_m(E)}\|\eta\|_{L^2_m(E)}$ and converges as $R\to\infty$. This proves the inequality in (i). Both integrals in \eqref{eq:weighted-ibp-R} converge, so $\ell:=\lim_{R\to\infty}m(R)\|\eta(R)\|_E^2$ exists. If $\ell>0$, then $m(R)\|\eta(R)\|_E^2\ge\ell/2$ for all large $R$, which contradicts $\int_0^\infty m\|\eta\|_E^2<\infty$. Hence $\ell=0$. This proves (ii) and the identity in (i).

Since $m$ is nonincreasing and integrable, $m(s)\to0$ and $\int_0^\infty[-m'(s)]\,ds=m(0)$. The Cauchy--Schwarz inequality gives the bound in (iii). Integration of $(m\eta)'$ on $[0,R]$ gives
\[
 \int_0^Rm(s)\eta_s(s)\,ds=m(R)\eta(R)+\int_0^R[-m'(s)]\eta(s)\,ds .
\]
We have $m(R)\|\eta(R)\|_E\le m(0)^{1/2}\bigl(m(R)\|\eta(R)\|_E^2\bigr)^{1/2}\to0$. Moreover $\int_0^\infty m\|\eta_s\|_E\le\kappa_0^{1/2}\|\eta_s\|_{L^2_m(E)}$, with $\kappa_0=\int_0^\infty m$. Letting $R\to\infty$ proves the identity in (iii).
\end{proof}

\begin{theorem}[Generation and dissipativity of the abstract coupled system]
\label{thm:variational-generation}
Under assumptions \eqref{eq:variational-basic-assumptions} and the
hypotheses on $m$ above, $\mathscr A_m$ with domain
\eqref{eq:variational-domain} is maximally dissipative on $\mathscr X$.
It generates a contraction $C_0$-semigroup. For classical solutions,
\begin{equation}
 \begin{aligned}
 \frac{d}{dt}\frac12\|(u(t),u_t(t),z(t))\|_{\mathscr X}^2
 &=-\Real d(u_t(t),u_t(t))\\
 &\quad-\frac\gamma2\int_0^\infty[-m'(s)]\|z(t,s)\|_Y^2\,ds.
 \end{aligned}
 \label{eq:variational-dissipation}
\end{equation}
\end{theorem}
\begin{proof}
Let $U=(u,v,z)\in D(\mathscr A_m)$.
The mechanical terms cancel in the real part of the energy pairing.
The two mixed terms $-\gamma(Jz,\ell(v))_Y$ and
$\gamma(\ell(v),Jz)_Y$ also cancel in real part. Hence
\[
 \Real\langle\mathscr A_mU,U\rangle_{\mathscr X}
 =-\Real d(v,v)-\gamma\Real\int_0^\infty m(s)(z_s(s),z(s))_Y\,ds
 \le0.
\]
Since $z\in D(\mathcal G)$, Lemma~\ref{lem:weighted-ibp}(i) with $E=Y$ gives the last term as
$-\gamma\int_0^\infty[-m'(s)]\|z(s)\|_Y^2ds/2$. This proves dissipativity. For a classical solution, $U\in C^1([0,\infty);\mathscr X)$ and $\frac{d}{dt}\frac12\|U\|_{\mathscr X}^2=\Real\langle\mathscr A_mU,U\rangle_{\mathscr X}$. This gives the
identity \eqref{eq:variational-dissipation}.

To show surjectivity of $\lambda I-\mathscr A_m$ for $\lambda>0$, let
$(f,g,k)\in\mathscr X$. The first equation gives $v=\lambda u-f$.
The history equation has the unique solution given by
\[
 z(s)=c_\lambda(s)\ell(v)+(R_\lambda k)(s),
\]
where
\[
 c_\lambda(s)=\frac{1-e^{-\lambda s}}\lambda\quad\text{and}\quad
 (R_\lambda k)(s)=\int_0^s e^{-\lambda(s-r)}k(r)\,dr .
\]
The contraction property of weighted translation implies
$\|R_\lambda k\|_{\mathcal M_m}\le\lambda^{-1}\|k\|_{\mathcal M_m}$.
Set $\mu_\lambda=\int_0^\infty m(s)c_\lambda(s)\,ds\ge0$ and
$k_\lambda=JR_\lambda k$. Substitution in the mechanical equation
leads to the following variational problem:
\begin{align*}
 \lambda^2(u,\psi)_H+a(u,\psi)+\lambda d(u,\psi)
 &+\gamma\lambda\mu_\lambda(\ell(u),\ell(\psi))_Y\\
 &= (g+\lambda f,\psi)_H+d(f,\psi)\\
 &\quad+\gamma(\mu_\lambda\ell(f)-k_\lambda,\ell(\psi))_Y.
\end{align*}
Its left-hand form is coercive on $V$ by
\eqref{eq:variational-basic-assumptions}. Lax--Milgram yields a
unique $u\in V$. Then $v\in V$ and $z\in D(\mathcal G)$.
The equation in $V'$ implies
$A_au+\mathcal D_dv+\gamma\ell^*Jz=g-\lambda v\in H$, so $U\in D(\mathscr A_m)$.
Consequently $\operatorname{Ran}(\lambda I-\mathscr A_m)=\mathscr X$.
It remains to check density of the domain. If
$F\in D(\mathscr A_m)^\perp$, surjectivity provides
$U\in D(\mathscr A_m)$ such that
$F=(\lambda I-\mathscr A_m)U$. Since $F$ is orthogonal to
$D(\mathscr A_m)$ and $U\in D(\mathscr A_m)$,
\[
 0=\Real\langle F,U\rangle_{\mathscr X}
 =\lambda\|U\|_{\mathscr X}^2
  -\Real\langle\mathscr A_mU,U\rangle_{\mathscr X}
 \ge\lambda\|U\|_{\mathscr X}^2.
\]
Thus $U=F=0$ and $D(\mathscr A_m)$ is dense in $\mathscr X$.
The Lumer--Phillips theorem therefore gives the asserted
contraction semigroup. The energy identity follows from the
calculation above for classical solutions.
\end{proof}

\subsection{Dissipation and multiplier hypotheses}
\label{sec:abstract-mechanical}

Along a trajectory $U(t)=(u(t),u_t(t),z(t))$ we use the memory energy, the memory dissipation, the mechanical energy and the total energy, defined by
\[
 M(t)=\int_0^\infty m(s)\|z(t,s)\|_Y^2\,ds\quad\text{and}\quad
 D(t)=\int_0^\infty[-m'(s)]\|z(t,s)\|_Y^2\,ds,
\]
\[
 E_{\rm mec}(t)=\tfrac12\bigl[a(u(t),u(t))+\|u_t(t)\|_H^2\bigr]\quad\text{and}\quad
 E(t)=\tfrac12\|U(t)\|_{\mathscr X}^2=E_{\rm mec}(t)+\tfrac{\gamma}{2}M(t).
\]
By \eqref{eq:variational-dissipation}, $E'(t)=-\Real d(u_t,u_t)-\frac{\gamma}{2}D(t)$ for classical solutions.

The decay argument combines the dissipation identity with a
mechanical multiplier estimate. The latter controls the part of the
mechanical energy that is coupled to the dissipative terms.

\begin{assumption}[Dissipation and mechanical multiplier estimates]
\label{ass:decay-transfer}
For each classical trajectory, set $E_0(t)=E_{\rm mec}(t)$ and
$\varphi(t)=\ell(u(t))$. The total energy is $E=E_0+\tfrac{\gamma}{2}M$, as above.
We require the following estimates, with constants independent of time
and of the initial state.

\smallskip
\noindent\textup{Energy dissipation.}
There exist $c_D>0$ and a nonnegative function $D_0(t)$ such that
\begin{equation}
\label{eq:H1-transfer}\tag{H1}
 E'(t)\le -D_0(t)-c_DD(t).
\end{equation}

\smallskip
\noindent\textup{Mechanical multiplier estimate.}
There exist a differentiable real-valued functional $\Psi(t)$ and
constants $C_D,C_M,C_\Psi>0$ such that
\begin{equation}
\label{eq:H2-transfer}\tag{H2}
\begin{aligned}
 \Psi'(t)&\le -E_0(t)+C_DD_0(t)+C_MM(t),\\
 |\Psi(t)|&\le C_\Psi E(t).
\end{aligned}
\end{equation}

\smallskip
\noindent\textup{Bounded mechanical observation.}
There exists $C_\varphi>0$ such that
\begin{equation}
\label{eq:H3-transfer}\tag{H3}
 \|\varphi(t)\|_Y^2=\|\ell(u(t))\|_Y^2\le C_\varphi E_0(t).
\end{equation}
\end{assumption}

The quantity $D_0$ may include a nonnegative portion of the memory dissipation.

For the variational system, \eqref{eq:variational-dissipation}
verifies (H1) with $D_0=\Real d(v,v)$ and $c_D=\gamma/2$.
When $d=0$, one may instead take $D_0=\gamma D/4$ and
$c_D=\gamma/4$. Coercivity of $a$ and boundedness of
$\ell:V\to Y$ yield (H3). It remains to establish (H2) under
suitable structural conditions, as in the next two propositions.

\begin{proposition}[One sufficient mechanical multiplier condition]
\label{prop:abstract-multiplier}
Assume, in addition, that for some $c_0,C_d>0$,
\begin{equation*}
 \Real d(v,v)\ge c_0\|v\|_H^2
 \quad\text{and}\quad
 |d(v,u)|\le C_d\bigl(\Real d(v,v)\bigr)^{1/2}a(u,u)^{1/2},
 \quad u,v\in V.
\end{equation*}
For a classical solution of \eqref{eq:variational-system}, set
$\Psi(t)=\Real(u(t),u_t(t))_H$.
There exist $c,C_1,C_2>0$, depending only on the structural constants,
such that
\begin{equation}
 \Psi'(t)\le-cE_{\rm mec}(t)
 +C_1\Real d(u_t(t),u_t(t))
 +C_2\int_0^\infty m(s)\|z(t,s)\|_Y^2\,ds.
 \label{eq:abstract-multiplier}
\end{equation}
\end{proposition}
\begin{proof}
Writing $v=u_t$ and testing the mechanical equation with $u$ gives
\[
 \Psi'=\|v\|_H^2-a(u,u)-\Real d(v,u)
 -\gamma\Real(Jz,\ell(u))_Y.
\]
Set $C_\ell=\|\ell\|_{\mathcal L(V,Y)}/\sqrt{c_a}$.
Then $\|\ell(u)\|_Y\le C_\ell a(u,u)^{1/2}$.
The weighted Cauchy--Schwarz inequality yields
$\|Jz\|_Y\le\kappa_0^{1/2}\|z\|_{\mathcal M_m}$.
Young's inequality gives
\begin{align*}
 |d(v,u)|&\le\tfrac14a(u,u)+C_d^2\Real d(v,v)\\
 \text{and}\quad\gamma\big|(Jz,\ell(u))_Y\big|
 &\le\tfrac14a(u,u)+\gamma^2\kappa_0 C_\ell^2M(t).
\end{align*}
Consequently,
\[
 \Psi'\le- E_{\rm mec}
 +\left(\frac{3}{2c_0}+C_d^2\right)\Real d(v,v)
 +\gamma^2\kappa_0 C_\ell^2M(t).
\]
The continuous embedding $V\hookrightarrow H$ also yields
$|\Psi|\le C E_{\rm mec}$. This proves
\eqref{eq:abstract-multiplier}, with the multiplier coefficients
expressed in terms of the displayed structural constants.
\end{proof}

Under an interior damping coercivity condition, normalize \eqref{eq:abstract-multiplier} by the positive coefficient $c$ and set $D_0=\Real d(v,v)$. This is exactly (H2) of Assumption~\ref{ass:decay-transfer}. The natural energy identity supplies (H1), while boundedness of $\ell$ on $V$ supplies (H3). Thus Proposition~\ref{prop:abstract-multiplier} verifies the decay hypotheses for oscillators, wave equations with interior friction and plates with rank-one memory and distributed damping. When no instantaneous damping is present, the mechanical energy must be recovered from the memory itself, as follows.

\begin{proposition}[Mechanical multiplier for a coercive memory observation]
\label{prop:full-memory-multiplier}
Assume $d=0$ and let $m$ be a positive, nonincreasing function in
$W^{1,1}(\R_+)$, with $m(0)<\infty$ and
$\kappa_0:=\int_0^\infty m(s)\,ds>0$. In particular,
$\int_0^\infty[-m'(s)]\,ds=m(0)$. Assume also the following condition:
\begin{equation}
 \|\ell u\|_Y\ge c_\ell\|u\|_V,\qquad u\in V,
 \label{eq:coercive-memory-observation}
\end{equation}
for some $c_\ell>0$, and assume $\operatorname{Ran}\ell=Y$, so that
$\ell:V\to Y$ is a bounded isomorphism. For every classical
solution of \eqref{eq:variational-system}, there exists a differentiable
real-valued functional $\widehat\Psi$ such that
\begin{equation}
 \widehat\Psi'(t)\le -E_{\rm mec}(t)+C_M M(t)+C_D D(t)
 \quad\text{and}\quad |\widehat\Psi(t)|\le C_\Psi E(t).
 \label{eq:full-memory-multiplier-estimate}
\end{equation}
\end{proposition}
\begin{proof}
Surjectivity and the lower bound
\eqref{eq:coercive-memory-observation} give a bounded inverse
$L:=\ell^{-1}:Y\to V$. For a classical history let $\eta(t,s)=Lz(t,s)$, and define
\[
 J\eta(t)=\int_0^\infty m(s)\eta(t,s)\,ds=LJz(t)\quad\text{and}\quad
 K(t)=\int_0^\infty[-m'(s)]\eta(t,s)\,ds.
\]
Let $U=(u,v,z)$ be a classical solution, so that $U\in C^1([0,\infty);\mathscr X)$ and $U(t)\in D(\mathscr A_m)$ for every $t\ge0$. Thus $u\in C^1([0,\infty);V)$ with $u_t=v$, and $v\in C^1([0,\infty);H)$ with
\[
 v_t=-A_au-\gamma\ell^*Jz\in H .
\]
Moreover $z\in C^1([0,\infty);\mathcal M_m)$ with $z_t=\ell(v)-z_s$, where $\ell(v)$ denotes a constant history. For each $t$, $z(t)\in D(\mathcal G)$. Since $L$ is bounded, $\eta(t)=Lz(t)$ satisfies the hypotheses of Lemma~\ref{lem:weighted-ibp} with $E=V$. Hence $K(t)$ is well defined in $V$ and $LJz_s(t)=K(t)$, by Lemma~\ref{lem:weighted-ibp}(iii). The map $J:\mathcal M_m\to Y$ is bounded. Therefore $J\eta=LJz\in C^1([0,\infty);V)$, and
\begin{equation}
 (J\eta)'(t)=LJ\bigl(\ell(v(t))\bigr)-LJz_s(t)=\kappa_0v(t)-K(t)\quad\text{with}\quad
 \kappa_0=\int_0^\infty m(s)\,ds.
 \label{eq:lifted-history-moment-derivative}
\end{equation}
These regularity properties justify the differentiation of the functionals below. They also justify the duality $\langle v_t,\psi\rangle_{V',V}=(v_t,\psi)_H$ for $\psi\in V$.

Define the two real-valued functionals by
\[
 \Psi_0(t)=\Real(u(t),v(t))_H\quad\text{and}\quad
 \Phi(t)=-\Real(v(t),J\eta(t))_H.
\]
Set
$C_{\ell,a}=\|\ell\|_{\mathcal L(V,Y)}c_a^{-1/2}$.
Testing the mechanical equation with $u$ and applying Young's
inequality, the continuity of $\ell$ and
$\|Jz\|_Y^2\le\kappa_0 M$, we obtain
\begin{equation}
 \Psi_0'(t)=\|v\|_H^2-a(u,u)
 -\gamma\Real(Jz,\ell u)_Y
 \le \|v\|_H^2-\tfrac12a(u,u)
 +\tfrac12\gamma^2\kappa_0C_{\ell,a}^2 M(t).
 \label{eq:full-memory-first-multiplier}
\end{equation}
On the other hand, testing the equation with $J\eta\in V$ and using
\eqref{eq:lifted-history-moment-derivative} gives the exact identity
\begin{equation*}
 \Phi'(t)=-\kappa_0\|v\|_H^2
 +\Real a(u,J\eta)+\gamma\|Jz\|_Y^2
 +\Real(v,K)_H.
\end{equation*}
We make the dependence of the constants explicit. Denote the
continuity bound of $a$ by $C_a$, the norm of the embedding
$V\hookrightarrow H$ by $C_{VH}$, and define the constants
\[
 C_L=\|L\|_{\mathcal L(Y,V)}\quad\text{and}\quad
 A=C_aC_L\kappa_0^{1/2}c_a^{-1/2}.
\]
For every $\varepsilon>0$,
\[
 |a(u,J\eta)|\le\varepsilon a(u,u)
              +\frac{A^2}{4\varepsilon} M(t)\quad\text{and}\quad
 \gamma\|Jz\|_Y^2\le\gamma\kappa_0M(t).
\]
Moreover,
\[
 \|K\|_H\le C_{VH}C_Lm(0)^{1/2}D(t)^{1/2}.
\]
Using Young's inequality with parameter $\kappa_0/2$ gives
\[
 |(v,K)_H|\le\frac{\kappa_0}{2}\|v\|_H^2+
   \frac{C_{VH}^2C_L^2m(0)}{2\kappa_0}D(t).
\]
Consequently,
\begin{equation}
 \Phi'(t)\le-\tfrac{\kappa_0}{2}\|v\|_H^2
 +\varepsilon a(u,u)
 +\left(\frac{A^2}{4\varepsilon}+\gamma\kappa_0\right)M(t)
 +\frac{C_{VH}^2C_L^2m(0)}{2\kappa_0}D(t).
 \label{eq:full-memory-second-multiplier}
\end{equation}
Set $N=4/\kappa_0$, choose $\varepsilon=1/(4N)$, and put
$\widehat\Psi=2(\Psi_0+N\Phi)$. Combining
\eqref{eq:full-memory-first-multiplier} and
\eqref{eq:full-memory-second-multiplier}, we find
\[
 \widehat\Psi'\le-2\|v\|_H^2-\tfrac12a(u,u)
       +C(M+D)\le-E_{\rm mec}+C(M+D).
\]
Finally $\|J\eta\|_H\le C M^{1/2}$, hence
\[
 |\Psi_0|\le C E_{\rm mec}\quad\text{and}\quad
 |\Phi|\le \|v\|_H\|J\eta\|_H\le C E.
\]
This proves \eqref{eq:full-memory-multiplier-estimate} with constants
independent of the trajectory. Since $d=0$,
the energy identity is $E'=-\gamma D/2$. Taking
$D_0=\gamma D/4$ and $c_D=\gamma/4$ verifies (H1), and the
term $C_DD_0$ in (H2) absorbs $CD$. Boundedness of $\ell$ on $V$
verifies (H3). This completes the proof.
\end{proof}

\begin{remark}[Memory dissipation in the multiplier hypothesis]
\label{rem:full-memory-h2}
The proof of Proposition~\ref{prop:full-memory-multiplier}
uses only the positivity, monotonicity and integrability of $m$.
The polynomial form of $m$ enters only through the decay rates.
The compensated multiplier involves the memory state and satisfies
$|\widehat\Psi|\le CE$, which is the bound used in
Lemma~\ref{lem:dyadic-step}.
\end{remark}

\subsection{Inherited and generated histories for Hilbert-valued observations}
\label{sec:abstract-inherited}

Let $\varphi\in C^1([0,\infty);Y)$ and let $z$ solve
$z_t+z_s=\varphi'(t)$ with $z(t,0)=0$ and initial history $z_0$.
The characteristic formula reads as follows:
\begin{equation}
 z(t,s)=
 \begin{cases}
 \varphi(t)-\varphi(t-s),&0\le s\le t,\\
 \varphi(t)+\xi_0(s-t),&s>t,
 \end{cases}
 \quad\text{with}\quad \xi_0(s):=z_0(s)-\varphi(0).
 \label{eq:inherited-history-formula}
\end{equation}
The centered history $\xi_0$ determines the inherited contribution for any initial observed displacement $\varphi(0)$.

We denote by $R(t)$ the right-translation semigroup on $Y$-valued histories, defined by
\[
 (R(t)z)(s)=\begin{cases}0,&0\le s\le t,\\ z(s-t),&s>t.\end{cases}
\]
For the density $m$, the inherited energy is defined by
\begin{equation*}
 H_{\xi_0}(t)=\int_0^\infty m(t+r)\|\xi_0(r)\|_Y^2\,dr=\|R(t)\xi_0\|^2_{\mathcal M_m}.
\end{equation*}

\subsection{The two-sided transfer principle}
\label{sec:general-transfer}

This subsection contains the main result. The decay from centered weighted histories is compared with the exact squared norm of the weighted right translation, for general densities and weights. In this subsection $m$ is a general density. It satisfies the hypotheses of Theorem~\ref{thm:variational-generation}, and $m(0)<\infty$. We put
\[
 \kappa(y)=\int_y^\infty m(s)\,ds,\quad\text{for all }y\ge0 .
\]

\begin{definition}
A nonincreasing function $f:\R_+\to(0,\infty)$ has \emph{decay of order at most $\lambda\ge0$} if there is $C_f\ge1$ such that
\begin{equation}
 f(s)\le C_f\Bigl(\frac ts\Bigr)^{\lambda}f(t),\quad\text{for all }1\le s\le t .
 \label{eq:decay-order}
\end{equation}
\end{definition}

A nonincreasing function $f$ has decay of order at most $\lambda$ for some $\lambda\ge0$ if and only if it satisfies the doubling property $f(s)\le Cf(2s)$ for $s\ge1$. Indeed, \eqref{eq:decay-order} with $t=2s$ gives doubling with $C=C_f2^{\lambda}$. Conversely, if $2^ks\le t<2^{k+1}s$, iteration and monotonicity give $f(s)\le C^{k+1}f(2^{k+1}s)\le C^{k+1}f(t)$, which is \eqref{eq:decay-order} with $\lambda=\log_2C$ and $C_f=C$. The exponent is not preserved. For instance, $f(t)=(1+t)^{-2}$ is doubling but has no decay of order at most $1$. Exponentially decaying functions are not doubling.

\begin{assumption}[Kernel and weight]
\label{ass:kernel-weight}
\begin{enumerate}[label=\textup{(K\arabic*)},leftmargin=*]
\item There is $a_*>0$ such that $-m'(s)\ge a_*m(s)/(1+s)$, for almost every $s\ge0$.
\item The density $m$ has decay of order at most $\lambda$, with constant $C_m$.
\item There is $K_m>0$ such that $\kappa(y)\le K_m(1+y)m(y)$, for all $y\ge0$.
\end{enumerate}
\begin{enumerate}[label=\textup{(W)},leftmargin=*]
\item The weight $w:\R_+\to[1,\infty)$ is measurable, $mw\in L^1_{\rm loc}(\R_+)$, and for every $T_*>0$
\[
 \Lambda_{T_*}:=\sup_{0\le t\le T_*}\ \operatorname*{ess\,sup}_{r\ge0}\frac{m(r+t)w(r+t)}{m(r)w(r)}<\infty .
\]
\end{enumerate}
\end{assumption}

Hypothesis (K3) is used only in the second part of Theorem~\ref{thm:general-transfer}. We define the spaces $\mathcal M_w$ and $\mathscr X_w$ by
\[
 \mathcal M_w=L^2(\R_+,m(s)w(s)\,ds;Y)\quad\text{and}\quad
 \mathscr X_w=\{(u,v,z)\in\mathscr X:z-\ell(u)\in\mathcal M_w\},
\]
with the norm given by $\|U\|_{\mathscr X_w}^2=\|U\|_{\mathscr X}^2+\|z-\ell(u)\|_{\mathcal M_w}^2$. Since $w\ge1$, the embedding $\mathcal M_w\hookrightarrow\mathcal M_m$ is continuous. Constant histories belong to $\mathcal M_m$ because $m$ is integrable. Hence a Cauchy sequence $U_n$ in $\mathscr X_w$ converges in $\mathscr X$ to some $U$, its centered histories converge in $\mathcal M_w$ to some $\xi$, and $\xi=z-\ell(u)$ in $\mathcal M_m$. Thus $\mathscr X_w$ is complete. The squared translation norm is defined by
\begin{equation}
 \rho_w(t)=\operatorname*{ess\,sup}_{r\ge0}\frac{m(r+t)}{m(r)w(r)},\quad\text{for all }t\ge0 .
 \label{eq:rho-w}
\end{equation}

\begin{proposition}[Exact norm of the weighted translation]
\label{prop:general-translation-norm}
For every $t\ge0$, $\|R(t)\|^2_{\mathcal L(\mathcal M_w,\mathcal M_m)}=\rho_w(t)$.
\end{proposition}
\begin{proof}
A change of variables gives
\[
 \|R(t)\xi\|_{\mathcal M_m}^2=\int_0^\infty\frac{m(r+t)}{m(r)w(r)}\,m(r)w(r)\|\xi(r)\|_Y^2\,dr\le\rho_w(t)\|\xi\|_{\mathcal M_w}^2 .
\]
Let $\varepsilon>0$. Since $mw$ is positive and locally integrable, there is a measurable set $A_\varepsilon$ of positive and finite $mw$-measure on which the quotient is at least $\rho_w(t)-\varepsilon$. Fix a unit vector $e\in Y$ and put
\[
 \xi_\varepsilon=\Bigl(\int_{A_\varepsilon}m(r)w(r)\,dr\Bigr)^{-1/2}\mathbf 1_{A_\varepsilon}e .
\]
Then $\|\xi_\varepsilon\|_{\mathcal M_w}=1$ and $\|R(t)\xi_\varepsilon\|^2_{\mathcal M_m}\ge\rho_w(t)-\varepsilon$.
\end{proof}

The decay estimates below concern the operator $T(t)$ from $\mathscr X_w$ into $\mathscr X$. The next proposition shows that $\mathscr X_w$ is invariant under the dynamics, so that the weighted spaces are evolution spaces of the semigroup.

\begin{proposition}[Invariance of $\mathscr X_w$]
\label{prop:general-invariance}
Under (W), the semigroup $T(t)$ leaves $\mathscr X_w$ invariant, and its restriction is strongly continuous on $\mathscr X_w$. In particular $D(\mathscr A_m)\cap\mathscr X_w$ is dense in $\mathscr X_w$.
\end{proposition}
\begin{proof}
Fix $U_0=(u_0,v_0,z_0)\in\mathscr X_w$, set $\xi_0=z_0-\ell(u_0)$, and write $T(t)U_0=(u(t),v(t),z(t))$ and $\varphi(t)=\ell(u(t))$. For classical solutions, \eqref{eq:inherited-history-formula} gives
\[
 z(t)-\varphi(t)=R(t)\xi_0-K_t\varphi\quad\text{with}\quad K_t\varphi(s)=\mathbf 1_{(0,t)}(s)\varphi(t-s).
\]
This identity extends to mild solutions by strong continuity of $T(t)$ in $\mathscr X$, continuity of $\ell:V\to Y$ and boundedness of $R(t)$ on $\mathcal M_m$. Let $T_*>0$ and $0\le t\le T_*$. By (W),
\[
 \|R(t)\xi_0\|^2_{\mathcal M_w}=\int_0^\infty m(r+t)w(r+t)\|\xi_0(r)\|_Y^2\,dr\le\Lambda_{T_*}\|\xi_0\|^2_{\mathcal M_w}.
\]
Coercivity of $a$, boundedness of $\ell$ and contractivity give $\sup_{r\ge0}\|\varphi(r)\|_Y^2\le C\|U_0\|^2_{\mathscr X}$. Hence
\[
 \|K_t\varphi\|^2_{\mathcal M_w}\le C\|U_0\|^2_{\mathscr X}\int_0^tm(s)w(s)\,ds .
\]
This proves invariance and local boundedness of the restricted semigroup. The space $C_c((0,\infty);Y)$ is dense in $\mathcal M_w$, because $mw$ is positive and locally integrable. Let $\xi\in C_c((0,\infty);Y)$ with support in $[a,b]$, $0<a<b$. For $0\le t<a$, the function $R(t)\xi$ is continuous with support in $[a,b+a]$, and $R(t)\xi\to\xi$ uniformly as $t\downarrow0$. Since $mw$ is integrable on $[0,b+a]$, this gives $R(t)\xi\to\xi$ in $\mathcal M_w$. The uniform bound on $R(t)$ extends this convergence to $\mathcal M_w$. Moreover $K_t\varphi\to0$ in $\mathcal M_w$. With the strong continuity of $T(t)$ in $\mathscr X$, this gives strong continuity in $\mathscr X_w$. The generator of the restriction has a dense domain in $\mathscr X_w$. Its difference quotients converge in $\mathscr X$, since $\mathscr X_w\hookrightarrow\mathscr X$ continuously. Hence its domain is contained in $D(\mathscr A_m)\cap\mathscr X_w$, which is therefore dense.
\end{proof}

\begin{lemma}[History splitting for a general density]
\label{lem:general-splitting}
Assume (K1). Let $z$ be given by \eqref{eq:inherited-history-formula}. Then, for every $t\ge0$ and $y\ge1$,
\begin{equation}
 M(t)\le\frac{2y}{a_*}D(t)+2H_{\xi_0}(t)+2\kappa(y)\|\varphi(t)\|_Y^2+2\int_y^tm(s)\|\varphi(t-s)\|_Y^2\,ds,
 \label{eq:general-splitting}
\end{equation}
where the last integral is zero if $y\ge t$. If $\sup_{r\ge0}\|\varphi(r)\|_Y^2\le B$, then
\begin{equation}
 M(t)\le\frac{2y}{a_*}D(t)+4B\kappa(y)+2H_{\xi_0}(t),\quad\text{for all }y\ge1 .
 \label{eq:general-splitting-first}
\end{equation}
Assume also (K3) and, for some $B>0$,
\[
 \|\varphi(t)\|_Y^2\le\frac{B}{1+t}\quad\text{and}\quad \int_0^\infty\|\varphi(r)\|_Y^2\,dr\le B .
\]
Then, for $1\le y\le1+t$,
\begin{equation}
 M(t)\le\frac{2y}{a_*}D(t)+(4K_m+2)Bm(y)+2H_{\xi_0}(t).
 \label{eq:general-splitting-second}
\end{equation}
\end{lemma}
\begin{proof}
By (K1), $m(s)\le(1+s)[-m'(s)]/a_*$. Hence
\[
 \int_0^ym(s)\|z(t,s)\|_Y^2\,ds\le\frac{1+y}{a_*}D(t)\le\frac{2y}{a_*}D(t).
\]
On $(y,\infty)$ we split the characteristic formula at $s=t$ and use $\|a+b\|^2_Y\le2\|a\|^2_Y+2\|b\|^2_Y$. This gives \eqref{eq:general-splitting}. Under the uniform bound, the last integral is at most $B\kappa(y)$, which gives \eqref{eq:general-splitting-first}. Under the second set of hypotheses, monotonicity of $m$ gives
\[
 \int_y^tm(s)\|\varphi(t-s)\|_Y^2\,ds\le m(y)\int_0^\infty\|\varphi(r)\|_Y^2\,dr\le Bm(y).
\]
For $1\le y\le1+t$, (K3) gives $\kappa(y)\|\varphi(t)\|_Y^2\le K_mBm(y)(1+y)/(1+t)\le2K_mBm(y)$. This proves \eqref{eq:general-splitting-second}.
\end{proof}

\begin{lemma}[Inherited decay order]
\label{lem:inherited-order}
Under (K2), the functions $\kappa$ and $\rho_w$ are nonincreasing and have decay of order at most $\lambda$, with constant $C_m$. A sum of two functions of decay of order at most $\lambda$ has the same property, with the larger constant.
\end{lemma}
\begin{proof}
For each $r$, the quotient $m(r+t)/(m(r)w(r))$ is nonincreasing in $t$. Hence $\rho_w$ is nonincreasing. Let $1\le s\le t$ and $r\ge0$. Then $r+s\ge1$ and $(r+t)/(r+s)\le t/s$. Applying \eqref{eq:decay-order} to $m$ at the points $r+s\le r+t$ gives
\[
 m(r+s)\le C_m\Bigl(\frac ts\Bigr)^{\lambda}m(r+t).
\]
Dividing by $m(r)w(r)$ and taking the essential supremum proves the claim for $\rho_w$. For $\kappa$, the substitution $x=s\sigma$ gives $\kappa(s)=s\int_1^\infty m(s\sigma)\,d\sigma$. Since $1\le s\sigma\le t\sigma$,
\[
 \kappa(s)\le C_m\Bigl(\frac ts\Bigr)^{\lambda}\frac st\,t\int_1^\infty m(t\sigma)\,d\sigma\le C_m\Bigl(\frac ts\Bigr)^{\lambda}\kappa(t).
\]
The statement on sums is immediate.
\end{proof}

\begin{lemma}[Comparison of the kernel with the translation norm]
\label{lem:kernel-rho}
Under (K2) and (W), there is $C_{m,w}>0$ such that $m(t)\le C_{m,w}\rho_w(t)$, for all $t\ge0$.
\end{lemma}
\begin{proof}
Since $w$ is finite almost everywhere, there are $W_1<\infty$ and a set $A\subset(0,1)$ of positive measure on which $w\le W_1$. For almost every $r\in A$, we have $m(r+t)\ge m(t+1)$ and $m(r)w(r)\le m(0)W_1$. Hence $\rho_w(t)\ge m(t+1)/(m(0)W_1)$. For $t\ge1$, \eqref{eq:decay-order} with the points $t\le t+1\le2t$ gives $m(t)\le C_m2^{\lambda}m(t+1)$. For $0\le t<1$, we have $m(t)\le m(0)$ and $m(t+1)\ge m(2)>0$. The claim follows with
\[
 C_{m,w}=m(0)W_1\max\Bigl\{C_m2^{\lambda},\frac{m(0)}{m(2)}\Bigr\}.
\]
\end{proof}

\begin{lemma}[Linear dyadic recursion]
\label{lem:linear-recursion}
Let $f$ be nonincreasing with decay of order at most $\lambda$ and constant $C_f$. Let $S_0\ge1$, $N\ge0$, $L\ge0$ and $\delta\in(0,1)$ with $\delta2^{\lambda+1}\le1$. Let $e:\R_+\to[0,N]$ be nonincreasing with
\begin{equation}
 e(2S)\le\delta e(S)+LNf(S),\quad\text{for all }S\ge S_0 .
 \label{eq:linear-recursion}
\end{equation}
Then $e(t)\le C_{\rm rec}Nf(t)$ for every $t\ge0$, with
\[
 C_{\rm rec}=2^{\lambda}C_f^{2}\Bigl(\frac1{f(S_0)}+2^{\lambda+1}L\Bigr)+\frac1{f(2S_0)} .
\]
\end{lemma}
\begin{proof}
Put $a_k=e(2^kS_0)$ and $f_k=f(2^kS_0)$. Iteration of \eqref{eq:linear-recursion} gives
\[
 a_k\le\delta^ka_0+LN\sum_{j=0}^{k-1}\delta^{k-1-j}f_j .
\]
For $0\le j\le k$, \eqref{eq:decay-order} gives $f_j\le C_f2^{\lambda(k-j)}f_k$. Since $\delta2^{\lambda}\le\frac12$,
\[
 \sum_{j=0}^{k-1}\delta^{k-1-j}f_j\le C_f2^{\lambda}f_k\sum_{i\ge0}(\delta2^{\lambda})^i\le2^{\lambda+1}C_ff_k .
\]
Moreover $a_0\le N$, $\delta^k\le2^{-\lambda k}$ and $f_0\le C_f2^{\lambda k}f_k$. Hence $\delta^ka_0\le C_fNf_k/f(S_0)$. Therefore $a_k\le C_f(1/f(S_0)+2^{\lambda+1}L)Nf_k$. For $t\ge S_0$, choose $k$ with $2^kS_0\le t<2^{k+1}S_0$. Then $e(t)\le a_k$ and $f_k\le C_f2^{\lambda}f(t)$. For $0\le t<2S_0$, we use $e(t)\le N\le Nf(t)/f(2S_0)$.
\end{proof}

\begin{lemma}[One dyadic step]
\label{lem:dyadic-step}
Assume \eqref{eq:H1-transfer}, \eqref{eq:H2-transfer} and (K1). Let $g:[1,\infty)\to(0,\infty)$ be nonincreasing with decay of order at most $\lambda$ and constant $C_g$. Let $U$ be a classical trajectory and $B\ge0$, and assume that
\begin{equation}
 M(t)\le\frac{2y}{a_*}D(t)+C_1Bg(y)+2H_{\xi_0}(t)
 \label{eq:dyadic-step-input}
\end{equation}
for all $t\ge1$ and all $1\le y\le t$. Then for every $\delta\in(0,1)$ there are $\varepsilon\in(0,1]$ and $S_0\ge1$, depending only on $\delta$, $a_*$, $\gamma$ and the constants of Assumption~\ref{ass:decay-transfer}, such that
\[
 E(2S)\le\delta E(S)+C_2\bigl(Bg(S)+H_{\xi_0}(S)\bigr),\quad\text{for all }S\ge S_0,
\]
with $C_2=2C_M'(C_1C_g\varepsilon^{-\lambda}+2)$ and $C_M'=1+\gamma/2+C_M$.
\end{lemma}
\begin{proof}
Fix $S\ge1$. Integration of \eqref{eq:H2-transfer} on $[S,2S]$ gives
\[
 \int_S^{2S}E_0\,dt\le|\Psi(S)|+|\Psi(2S)|+C_D\int_S^{2S}D_0\,dt+C_M\int_S^{2S}M\,dt .
\]
The energy is nonincreasing. Hence $|\Psi(S)|+|\Psi(2S)|\le2C_\Psi E(S)$. Integration of \eqref{eq:H1-transfer} gives
\[
 \int_S^{2S}D_0\,dt\le E(S)\quad\text{and}\quad\int_S^{2S}D\,dt\le\frac{E(S)}{c_D}.
\]
Since $E=E_0+\frac\gamma2M$, we obtain
\[
 \int_S^{2S}E\,dt\le C_3E(S)+C_M'\int_S^{2S}M\,dt\quad\text{with}\quad C_3=2C_\Psi+C_D .
\]
Let $\varepsilon\in(0,1]$, $S\ge1/\varepsilon$ and $y=\varepsilon S$. Then $1\le y\le S\le t$ for $t\in[S,2S]$, so \eqref{eq:dyadic-step-input} holds on this interval. The function $H_{\xi_0}$ is nonincreasing. Integration of \eqref{eq:dyadic-step-input} gives
\[
 \int_S^{2S}M\,dt\le\frac{2\varepsilon S}{a_*c_D}E(S)+C_1BSg(\varepsilon S)+2SH_{\xi_0}(S).
\]
Since $SE(2S)\le\int_S^{2S}E\,dt$, division by $S$ yields
\[
 E(2S)\le\Bigl(\frac{C_3}{S}+\frac{2C_M'\varepsilon}{a_*c_D}\Bigr)E(S)+C_M'C_1Bg(\varepsilon S)+2C_M'H_{\xi_0}(S).
\]
Since $1\le\varepsilon S\le S$, we have $g(\varepsilon S)\le C_g\varepsilon^{-\lambda}g(S)$. The choice
\[
 \varepsilon=\min\Bigl\{1,\frac{\delta a_*c_D}{4C_M'}\Bigr\}\quad\text{and}\quad S_0=\max\Bigl\{\frac1\varepsilon,\frac{2C_3}{\delta}\Bigr\}
\]
makes the coefficient of $E(S)$ at most $\delta$ for $S\ge S_0$.
\end{proof}

\begin{lemma}[Kernel and tail]
\label{lem:kernel-tail}
Under (K2), there is $C_{m\kappa}>0$ such that $m(t)\le C_{m\kappa}\kappa(t)$, for all $t\ge0$.
\end{lemma}
\begin{proof}
For $t\ge1$, monotonicity and \eqref{eq:decay-order} give $\kappa(t)\ge\int_t^{2t}m\ge t\,m(2t)\ge m(t)/(C_m2^{\lambda})$. For $0\le t<1$, we have $m(t)\le m(0)$ and $\kappa(t)\ge\kappa(1)>0$. The claim holds with $C_{m\kappa}=\max\{C_m2^{\lambda},m(0)/\kappa(1)\}$.
\end{proof}

\begin{lemma}[Lower bound by the inherited and generated histories]
\label{lem:historical-lower}
Let $m$ satisfy the hypotheses of Theorem~\ref{thm:variational-generation}, and put $\kappa_0=\int_0^\infty m$. For every $U_0\in\mathscr X$ and every $t\ge0$,
\begin{equation}
 E(t)\ge c_*\left(H_{\xi_0}(t)+m(t)\int_0^t\|\varphi(r)\|_Y^2\,dr\right)
 \quad\text{with}\quad c_*=\Bigl(\frac8\gamma+4\kappa_0C_\varphi\Bigr)^{-1}.
 \label{eq:sharp-historical-lower}
\end{equation}
\end{lemma}
\begin{proof}
The characteristic formula \eqref{eq:inherited-history-formula} holds for mild solutions. For $s>t$ it gives $\xi_0(s-t)=z(t,s)-\varphi(t)$. Hence
\[
 H_{\xi_0}(t)=\int_t^\infty m(s)\|z(t,s)-\varphi(t)\|_Y^2\,ds\le2M(t)+2\kappa(t)\|\varphi(t)\|_Y^2 .
\]
For $0<s\le t$ it gives $\varphi(t-s)=\varphi(t)-z(t,s)$. Since $m(s)\ge m(t)$ for $s\le t$,
\[
 m(t)\int_0^t\|\varphi(r)\|_Y^2\,dr\le\int_0^tm(s)\|\varphi(t-s)\|_Y^2\,ds\le2M(t)+2\bigl(\kappa_0-\kappa(t)\bigr)\|\varphi(t)\|_Y^2 .
\]
The sum of the two left-hand sides is at most $4M(t)+4\kappa_0\|\varphi(t)\|_Y^2$. Since $M\le2E/\gamma$ and $\|\varphi(t)\|_Y^2\le C_\varphi E_0(t)\le C_\varphi E(t)$ by \eqref{eq:H3-transfer}, this sum is at most $c_*^{-1}E(t)$.
\end{proof}

\begin{lemma}[Resolvent approximation and the inherited history source]
\label{lem:resolvent-history-source}
Let $m$ satisfy the hypotheses of Theorem~\ref{thm:variational-generation}, and let $\mathscr A_m$ be the generator of the contraction semigroup. For $U_0=(u_0,v_0,z_0)\in\mathscr X$,
put $\xi_0=z_0-\ell(u_0)$ and, for $\lambda\ge1$, define
\[
 U_\lambda=(u_\lambda,v_\lambda,z_\lambda)
 =\lambda(\lambda I-\mathscr A_m)^{-1}U_0
 \quad\text{with}\quad \xi_\lambda=z_\lambda-\ell(u_\lambda).
\]
Then $U_\lambda\in D(\mathscr A_m)$,
$\|U_\lambda\|_{\mathscr X}\le\|U_0\|_{\mathscr X}$ and
$U_\lambda\to U_0$ in $\mathscr X$ as $\lambda\to\infty$.
Moreover, for every $t\ge0$,
\begin{equation}
 H_{\xi_\lambda}(t)
 \le 2H_{\xi_0}(t)
 +\frac{C}{\lambda}\|U_0\|_{\mathscr X}^{2}
     m(t),
 \label{eq:resolvent-source-stability}
\end{equation}
where $C>0$ depends only on the fixed system coefficients and the norm of $\ell:V\to Y$.
\end{lemma}
\begin{proof}
Since $\mathscr A_m$ is maximally dissipative, the standard
resolvent estimate and the Yosida approximation give
\[
 \|U_\lambda\|_{\mathscr X}\le\|U_0\|_{\mathscr X}
 \quad\text{and}\quad U_\lambda\longrightarrow U_0\quad\hbox{in }\mathscr X.
\]
The first and third components of
$(\lambda I-\mathscr A_m)U_\lambda=\lambda U_0$ read
\[
 v_\lambda=\lambda(u_\lambda-u_0)\quad\text{and}\quad
 \partial_s z_\lambda+\lambda z_\lambda
 =\ell(v_\lambda)+\lambda z_0\quad\text{and}\quad z_\lambda(0)=0.
\]
Define the resolvent of the right-translation generator by
\[
 (R_\lambda f)(s)=\int_0^s e^{-\lambda(s-r)}f(r)\,dr.
\]
Since $z_0=\xi_0+\ell(u_0)$, direct integration gives
\[
 \begin{aligned}
 z_\lambda(s)
 &=\frac{1-e^{-\lambda s}}{\lambda}\ell(v_\lambda)
     +\lambda R_\lambda z_0(s)\\
 &=(1-e^{-\lambda s})\ell(u_\lambda)
     +\lambda R_\lambda\xi_0(s),
 \end{aligned}
\]
and consequently the exact identity
\begin{equation}
 \xi_\lambda(s)=\lambda R_\lambda\xi_0(s)
       -e^{-\lambda s}\ell(u_\lambda).
 \label{eq:resolvent-centered-history}
\end{equation}
The first term is an average of right translates:
\[
 \lambda R_\lambda\xi_0(s)
 =\int_0^s\lambda e^{-\lambda\sigma}
             \xi_0(s-\sigma)\,d\sigma.
\]
As $\int_0^s\lambda e^{-\lambda\sigma}\,d\sigma\le1$,
Jensen's inequality (with the missing mass placed at zero) yields
\[
 \|\lambda R_\lambda\xi_0(s)\|_Y^2
 \le \int_0^s\lambda e^{-\lambda\sigma}
                  \|\xi_0(s-\sigma)\|_Y^2\,d\sigma.
\]
Tonelli's theorem and the monotonicity of $m$ now give, for every
$t\ge0$,
\[
 \begin{aligned}
 H_{\lambda R_\lambda\xi_0}(t)
 &\le\int_0^\infty\lambda e^{-\lambda\sigma}
       \int_0^\infty m(t+\sigma+r)
                       \|\xi_0(r)\|_Y^2\,dr\,d\sigma\\
 &\le\left(\int_0^\infty\lambda e^{-\lambda\sigma}
       \,d\sigma\right)H_{\xi_0}(t)
 =H_{\xi_0}(t).
 \end{aligned}
\]
On the other hand, by boundedness of $\ell$, coercivity of $a$,
and the resolvent contraction estimate,
$\|\ell(u_\lambda)\|_Y^2\le C\|U_0\|_{\mathscr X}^2$.
Since $m$ is decreasing,
\[
 \int_0^\infty m(t+s)e^{-2\lambda s}\,ds
 \le\frac{m(t)}{2\lambda}.
\]
Combining these estimates with
$\|a+b\|_Y^2\le2\|a\|_Y^2+2\|b\|_Y^2$
in \eqref{eq:resolvent-centered-history} proves
\eqref{eq:resolvent-source-stability}.
\end{proof}

The next theorem is the core of the method. The inherited energy enters only through a prescribed upper envelope $\sigma$.

\begin{theorem}[Transfer from a prescribed inherited source]
\label{thm:source-transfer}
Assume the hypotheses of Theorem~\ref{thm:variational-generation}, Assumption~\ref{ass:decay-transfer}, (K1) and (K2). Let $\sigma:\R_+\to(0,\infty)$ be nonincreasing with decay of order at most $\lambda$. Let $U_0\in\mathscr X$ and $K\ge0$ satisfy
\begin{equation}
 H_{\xi_0}(t)\le K\sigma(t),\quad\text{for all }t\ge0 .
 \label{eq:source-envelope}
\end{equation}
Put $N=E(0)+K$. There is $C>0$, independent of $t$, $U_0$ and $K$, such that
\begin{equation}
 E(t)\le CN\bigl(\kappa(t)+\sigma(t)\bigr),\quad\text{for all }t\ge0 .
 \label{eq:source-first}
\end{equation}
Assume moreover (K3), $\int_0^\infty s\,m(s)\,ds<\infty$ and $\sigma\in L^1(\R_+)$. Then there is $C'>0$, independent of $t$, $U_0$ and $K$, such that
\begin{equation}
 E(t)\le C'N\bigl(m(t)+\sigma(t)\bigr),\quad\text{for all }t\ge0 .
 \label{eq:source-second}
\end{equation}
In all cases, $E(t)\ge c_*H_{\xi_0}(t)$ for $t\ge0$.
\end{theorem}
\begin{proof}
The lower bound is \eqref{eq:sharp-historical-lower} of Lemma~\ref{lem:historical-lower}. If $N=0$, then $U_0=0$. Assume $N>0$.

\emph{Classical data.} Let $U_0\in D(\mathscr A_m)$. Contractivity and \eqref{eq:H3-transfer} give $E(t)\le N$ and $\|\varphi(t)\|_Y^2\le C_\varphi N$. Fix $\delta=2^{-\lambda-1}$.

\emph{First pass.} Estimate \eqref{eq:general-splitting-first} with $B=C_\varphi N$ gives \eqref{eq:dyadic-step-input} with $g=\kappa$ and $C_1=4$. By Lemma~\ref{lem:inherited-order}, $\kappa$ has decay of order at most $\lambda$. Lemma~\ref{lem:dyadic-step} and \eqref{eq:source-envelope} yield
\[
 E(2S)\le\delta E(S)+C_2(C_\varphi+1)N\bigl(\kappa(S)+\sigma(S)\bigr),\quad\text{for all }S\ge S_0 .
\]
By Lemma~\ref{lem:inherited-order}, $f_1=\kappa+\sigma$ is nonincreasing with decay of order at most $\lambda$. Lemma~\ref{lem:linear-recursion} gives \eqref{eq:source-first}.

\emph{Second pass.} The function $f_1$ is nonincreasing and integrable, since $\int_0^\infty\kappa=\int_0^\infty s\,m(s)\,ds$. For $t>0$, we have $\frac t2f_1(t)\le\int_{t/2}^tf_1\le\|f_1\|_{L^1}$. Hence $(1+t)f_1(t)\le f_1(0)+2\|f_1\|_{L^1}$. By \eqref{eq:H3-transfer} and \eqref{eq:source-first},
\[
 \|\varphi(t)\|_Y^2\le\frac{B_2}{1+t}\quad\text{and}\quad\int_0^\infty\|\varphi(r)\|_Y^2\,dr\le B_2\quad\text{with}\quad B_2=C_\varphi CN\bigl(f_1(0)+2\|f_1\|_{L^1}\bigr).
\]
Estimate \eqref{eq:general-splitting-second} holds for $1\le y\le t$. It gives \eqref{eq:dyadic-step-input} with $g=m$, $B=B_2$ and $C_1=4K_m+2$. Lemma~\ref{lem:dyadic-step} yields
\[
 E(2S)\le\delta E(S)+C_4N\bigl(m(S)+\sigma(S)\bigr),\quad\text{for all }S\ge S_0 .
\]
Lemma~\ref{lem:linear-recursion} with $f=m+\sigma$ proves \eqref{eq:source-second}.

\emph{General data.} Let $U_0\in\mathscr X$ and $\lambda\ge1$, and let $U_\lambda$ be as in Lemma~\ref{lem:resolvent-history-source}. Then $U_\lambda\in D(\mathscr A_m)$, $E_\lambda(0)\le E(0)$, and
\[
 H_{\xi_\lambda}(t)\le2K\sigma(t)+C\|U_0\|^2_{\mathscr X}m(t)\le K_*\bigl(\sigma(t)+m(t)\bigr)\quad\text{with}\quad K_*=2K+2C\,E(0).
\]
The envelope $\sigma+m$ is nonincreasing, has decay of order at most $\lambda$, and is integrable when $\sigma$ is. The classical case gives, with constants independent of $\lambda$,
\[
 E_\lambda(t)\le CN\bigl(\kappa(t)+\sigma(t)+m(t)\bigr)\le C(1+C_{m\kappa})N\bigl(\kappa(t)+\sigma(t)\bigr),
\]
by Lemma~\ref{lem:kernel-tail}, and $E_\lambda(t)\le C'N(m(t)+\sigma(t))$ under the additional hypotheses. Here $E_\lambda$ is the energy of $T(t)U_\lambda$. The resolvent commutes with $T(t)$, so $T(t)U_\lambda\to T(t)U_0$ in $\mathscr X$ for every $t$. Letting $\lambda\to\infty$ proves \eqref{eq:source-first} and \eqref{eq:source-second}.
\end{proof}

For $U_0\in\mathscr X_w$ we set $N=E(0)+\|\xi_0\|^2_{\mathcal M_w}$. Then $N\le\|U_0\|^2_{\mathscr X_w}$.

\begin{theorem}[General weighted-history transfer]
\label{thm:general-transfer}
Assume the hypotheses of Theorem~\ref{thm:variational-generation}, Assumption~\ref{ass:decay-transfer}, (K1), (K2) and (W). There is $C>0$, independent of $t$ and $U_0$, such that every $U_0\in\mathscr X_w$ satisfies
\begin{equation}
 E(t)\le CN\bigl(\kappa(t)+\rho_w(t)\bigr),\quad\text{for all }t\ge0 .
 \label{eq:general-first}
\end{equation}
Assume moreover (K3), $\int_0^\infty s\,m(s)\,ds<\infty$ and $\rho_w\in L^1(\R_+)$. Then there is $C'>0$, independent of $t$ and $U_0$, such that
\begin{equation}
 E(t)\le C'N\rho_w(t),\quad\text{for all }t\ge0 .
 \label{eq:general-second}
\end{equation}
\end{theorem}
\begin{proof}
By Proposition~\ref{prop:general-translation-norm}, $H_{\xi_0}(t)=\|R(t)\xi_0\|^2_{\mathcal M_m}\le\|\xi_0\|^2_{\mathcal M_w}\rho_w(t)$. By Lemma~\ref{lem:inherited-order}, $\rho_w$ is nonincreasing with decay of order at most $\lambda$. Theorem~\ref{thm:source-transfer} with $\sigma=\rho_w$ and $K=\|\xi_0\|^2_{\mathcal M_w}$ gives \eqref{eq:general-first}. Under the additional hypotheses it gives $E(t)\le C'N(m(t)+\rho_w(t))$. Lemma~\ref{lem:kernel-rho} gives $m\le C_{m,w}\rho_w$, which proves \eqref{eq:general-second}.
\end{proof}

\begin{corollary}[Two-sided transfer]
\label{cor:general-two-sided}
Assume the hypotheses of the first part of Theorem~\ref{thm:general-transfer}. Assume that one of the following conditions holds.
\begin{enumerate}[label=\textup{(\roman*)},leftmargin=*]
\item (K3) holds, $\int_0^\infty s\,m(s)\,ds<\infty$ and $\rho_w\in L^1(\R_+)$.
\item There is $C_\kappa>0$ such that $\kappa(t)\le C_\kappa\rho_w(t)$, for all $t\ge0$.
\end{enumerate}
Then there are $0<c\le C$ such that
\[
 c\,\rho_w(t)\le\|T(t)\|^2_{\mathcal L(\mathscr X_w,\mathscr X)}\le C\rho_w(t),\quad\text{for all }t\ge0 .
\]
\end{corollary}
\begin{proof}
\emph{Lower bound.} By \eqref{eq:inherited-history-formula}, $\xi_0(s-t)=z(t,s)-\varphi(t)$ for $s>t$. Hence
\[
 H_{\xi_0}(t)=\int_t^\infty m(s)\|z(t,s)-\varphi(t)\|_Y^2\,ds\le2M(t)+2\kappa(t)\|\varphi(t)\|_Y^2 .
\]
Since $M\le2E/\gamma$, $\kappa(t)\le\kappa(0)$ and \eqref{eq:H3-transfer} holds, we obtain
\[
 H_{\xi_0}(t)\le C_HE(t)\quad\text{with}\quad C_H=\frac4\gamma+2\kappa(0)C_\varphi .
\]
The characteristic formula holds for mild solutions, so this bound holds for every $U_0\in\mathscr X$.
For $U_0=(0,0,\xi_0)$, we have $\|U_0\|^2_{\mathscr X_w}\le(\gamma+1)\|\xi_0\|^2_{\mathcal M_w}$ and $\|T(t)U_0\|^2_{\mathscr X}=2E(t)$. Taking the supremum over $\xi_0$ and applying Proposition~\ref{prop:general-translation-norm} gives the lower bound with $c=2/((\gamma+1)C_H)$.

\emph{Upper bound.} We have $\|T(t)U_0\|^2_{\mathscr X}=2E(t)$ and $N\le\|U_0\|^2_{\mathscr X_w}$. Under (i), \eqref{eq:general-second} applies. Under (ii), \eqref{eq:general-first} gives $E(t)\le C(1+C_\kappa)N\rho_w(t)$.
\end{proof}

The weight built from the tail of the kernel gives a universal upper bound for the translation norm. Jin, Liang and Xiao \cite[Theorem~3.1]{JinLiangXiao2019} proved the rate $t^{-1}$ for every nonincreasing integrable kernel, for coupled systems and data of higher regularity whose histories are weighted by $\kappa$. In the present framework, under (K3) and a finite first moment, this rate is an exact operator-norm law without additional regularity.

\begin{corollary}[The tail weight]
\label{cor:tail-weight}
Let $w=1+\kappa/m$. For every nonincreasing positive integrable $m$,
\[
 \rho_w(t)\le\min\Bigl\{1,\frac1t\Bigr\}\le\frac2{1+t},\quad\text{for all }t>0 .
\]
Assume the hypotheses of the first part of Theorem~\ref{thm:general-transfer} except (W), and assume (K3) and $\int_0^\infty s\,m(s)\,ds<\infty$. Then (W) holds, and there are $0<c\le C$ such that
\[
 \frac{c}{1+t}\le\|T(t)\|^2_{\mathcal L(\mathscr X_w,\mathscr X)}\le\frac{C}{1+t},\quad\text{for all }t\ge0 .
\]
\end{corollary}
\begin{proof}
We have $m(r)w(r)=m(r)+\kappa(r)\ge\kappa(r)\ge\int_r^{r+t}m\ge t\,m(r+t)$. Hence the quotient in \eqref{eq:rho-w} is at most $1/t$. It is also at most $1$, since $w\ge1$ and $m$ is nonincreasing. This proves the first claim.

The function $mw=m+\kappa$ is nonincreasing and locally integrable. Hence (W) holds with $\Lambda_{T_*}\le1$. For $t\ge1$, the choice $r=t$ in \eqref{eq:rho-w}, (K2) and (K3) give
\[
 \rho_w(t)\ge\frac{m(2t)}{m(t)+\kappa(t)}\ge\frac{1}{C_m2^{\lambda}\bigl(1+K_m(1+t)\bigr)}\ge\frac{1}{C_m2^{\lambda}(1+K_m)(1+t)} .
\]
For $0\le t<1$, monotonicity gives $\rho_w(t)\ge\rho_w(1)$. Moreover $\kappa$ is nonincreasing and integrable. Hence $(1+t)\kappa(t)\le\kappa(0)+2\|\kappa\|_{L^1}$. Thus $\kappa\le C_\kappa\rho_w$, and alternative (ii) of Corollary~\ref{cor:general-two-sided} applies.
\end{proof}

\begin{remark}[Relation with trajectory lower bounds]
\label{rem:conti-pata-lower}
Lemma~\ref{lem:kernel-rho} is a uniform counterpart of Proposition~1.1 of Conti and Pata \cite{ContiPata2020}. That proposition gives $\mu(t)\le C(u_0,v_0)E(t)$ for each nonzero solution of the Volterra equation without prehistory. Here $m(t)\le C_{m,w}\rho_w(t)$ shows that the saturated rate is contained in the translation norm. Hypothesis (K1) excludes kernels with flat zones or jumps, which are covered by \cite[Theorem~4.4]{ContiPata2020}. Whether the two-sided law persists without (K1) is open.
\end{remark}

\subsection{Polynomial kernels: translation estimates}\label{sec:historical-decay}
\label{sec:abstract-weighted-history}

We now specialize the transfer principle to polynomial kernels and weights. This case yields explicit exponents. It also admits an extension to prescribed inherited tails. We begin with the norm of the right translation between polynomially weighted memory spaces. Section~\ref{sec:abstract-decay-theorem} transfers it to the coupled system under Assumption~\ref{ass:decay-transfer}.

In the remainder of Section~\ref{sec:variational}, $m(s)=pb(1+s)^{-p-1}$ with $p>1$ and $b>0$, and denote the fixed kernel exponent by $q_*=p+1$. The right translation $R(t)$ is that of Section~\ref{sec:abstract-inherited}.
For $\theta>0$, define
\begin{equation*}
 \mathcal M_\theta=L^2\!\bigl(\R_+,m(s)(1+s)^{2\theta}\,ds;Y\bigr).
\end{equation*}
Since $(1+s)^{2\theta}\ge1$, one has $\mathcal M_\theta\hookrightarrow\mathcal M_m$ with $\|z\|_{\mathcal M_m}\le\|z\|_{\mathcal M_\theta}$. A nonzero constant history belongs to $\mathcal M_\theta$ only when $2\theta<p$. For this reason the weighted spaces below are imposed on the centered history $z-\ell(u)$ rather than on $z$.

\begin{proposition}[Sharp decay of weighted histories under translation]
\label{prop:abstract-weighted-transport}
For the right-translation semigroup $R(t)$ and every $t\ge0$,
\begin{equation}
 \|R(t)\|_{\mathcal L(\mathcal M_\theta,\mathcal M_m)}^2
 =\sup_{s\ge0}\frac{(1+s)^{q_*-2\theta}}{(1+s+t)^{q_*}}.
 \label{eq:weighted-translation-exact}
\end{equation}
Consequently there are $0<c_{p,\theta}\le C_{p,\theta}<\infty$ such that
\begin{gather*}
 c_{p,\theta}(1+t)^{-\alpha(p,\theta)}
 \le \|R(t)\|_{\mathcal L(\mathcal M_\theta,\mathcal M_m)}^2
 \le C_{p,\theta}(1+t)^{-\alpha(p,\theta)},\\
 \text{with}\quad \alpha(p,\theta)=\min\{2\theta,p+1\}.
\end{gather*}
In particular, the saturation threshold for the translation component is
$\theta_c=(p+1)/2$.
\end{proposition}

\begin{proof}
Changing variables in the translation norm gives
\[
 \|R(t)z\|_{\mathcal M_m}^2=\int_0^\infty m(s+t)\|z(s)\|_Y^2\,ds.
\]
The norm of this multiplication-weight embedding is the essential supremum of
$m(s+t)/[m(s)(1+s)^{2\theta}]$, which is the quantity in
\eqref{eq:weighted-translation-exact}. Continuity gives the displayed supremum.
Writing $x=1+s\ge1$, we maximize $x^{q_*-2\theta}/(x+t)^{q_*}$.
If $2\theta\ge q_*$, the expression is nonincreasing and its supremum is
$(1+t)^{-q_*}$. If $0<2\theta<q_*$, its unrestricted maximum occurs at
$x=(q_*-2\theta)t/(2\theta)$. For all sufficiently large $t$ this point belongs to $[1,\infty)$ and the maximum is a positive constant times $t^{-2\theta}$. Enlarging the comparison constants covers bounded $t$.
\end{proof}

\subsection{Polynomial kernels: the coupled system and prescribed histories}
\label{sec:abstract-decay-theorem}

We now transfer the localization of the initial history to the coupled energy. Throughout this subsection, $T(t)=e^{t\mathscr A_m}$ is the contraction semigroup of Theorem~\ref{thm:variational-generation} on $\mathscr X$. For $U(t)=T(t)U_0=(u(t),v(t),z(t))$, the observation is $\varphi(t)=\ell(u(t))$, the quantities $M$, $D$ and $E$ are those of Section~\ref{sec:abstract-mechanical}, and $H_{\xi_0}$ is defined in Section~\ref{sec:abstract-inherited} with $\xi_0=z_0-\ell(u_0)$. We also put
\[
 \kappa(s)=\int_s^\infty m(r)\,dr=b(1+s)^{-p}.
\]
Unless otherwise stated, the positive constants may depend on the fixed coefficients, on $p$ and on $\theta$, but not on time or on the initial state.

Define the weighted initial-data space, centered at the initial value
of the mechanical observation,
\begin{equation*}
 \mathscr X_\theta
 =\left\{(u,v,z)\in\mathscr X:
 z-\ell(u)\in\mathcal M_\theta\right\}.
\end{equation*}
It carries the norm defined by
\begin{equation*}
 \|U\|_{\mathscr X_\theta}^2
 =\|U\|_{\mathscr X}^2+\|z-\ell(u)\|_{\mathcal M_\theta}^2.
\end{equation*}
The centered variable $z-\ell(u)$ isolates the inherited history
and determines the decay exponent in the weighted norm.
The space $\mathscr X_\theta$ is complete for this norm:
if $U_n=(u_n,v_n,z_n)$ is Cauchy in $\mathscr X_\theta$,
then $U_n\to U$ in $\mathscr X$ and
$z_n-\ell(u_n)\to\xi$ in $\mathcal M_\theta$.
Since $\mathcal M_\theta\hookrightarrow\mathcal M_m$ continuously
and $u\mapsto\ell(u)$ defines a constant $Y$-valued history
in $\mathcal M_m$, we have $z-\ell(u)=\xi$ in
$\mathcal M_m$. Hence $U\in\mathscr X_\theta$ and the
convergence holds in the full weighted norm.

The polynomial weight satisfies (W), with $\Lambda_{T_*}\le(1+T_*)^{\max\{2\theta-p-1,0\}}$. By Proposition~\ref{prop:general-invariance}, the semigroup $T(t)$ leaves $\mathscr X_\theta$ invariant, its restriction is strongly continuous on $\mathscr X_\theta$, and $D(\mathscr A_m)\cap\mathscr X_\theta$ is dense in $\mathscr X_\theta$.

\begin{theorem}[Sharp abstract transfer and saturation threshold]
\label{thm:sharp-abstract-transfer}
Let the contractive coupled system satisfy Assumption~\ref{ass:decay-transfer}, with
$m(s)=pb(1+s)^{-p-1}$, $p>1$, $b>0$. Then for every $\theta>0$,
\begin{equation}
 E(t)\le C_{p,\theta}(1+t)^{-\alpha(p,\theta)}
 \bigl(E(0)+\|\xi_0\|_{\mathcal M_\theta}^2\bigr)
 \quad\text{with}\quad \alpha(p,\theta)=\min\{2\theta,p+1\},
 \label{eq:sharp-abstract-transfer-upper}
\end{equation}
for every $U_0\in\mathscr X_\theta$, with $\xi_0=z_0-\ell(u_0)$.
Moreover, with the energy norm induced by $E$, there exist
$0<c_{p,\theta}\le C_{p,\theta}$ such that
\begin{equation}
 c_{p,\theta}(1+t)^{-\alpha(p,\theta)}
 \le\|T(t)\|_{\mathcal L(\mathscr X_\theta,\mathscr X)}^2
 \le C_{p,\theta}(1+t)^{-\alpha(p,\theta)},\quad\text{for all }t\ge0.
 \label{eq:sharp-abstract-transfer-bilateral}
\end{equation}
In particular, the exponent $\alpha(p,\theta)$ cannot be improved.
\end{theorem}
\begin{proof}
We verify the hypotheses of Corollary~\ref{cor:general-two-sided} for $m(s)=pb(1+s)^{-p-1}$ and $w(s)=(1+s)^{2\theta}$. Hypothesis (K1) holds with $a_*=p+1$, since $-m'(s)=(p+1)m(s)/(1+s)$. Hypothesis (K2) holds with $\lambda=p+1$ and $C_m=1$, because $(1+t)/(1+s)\le t/s$ for $1\le s\le t$. Hypothesis (K3) holds with $K_m=1/p$, since $\kappa(y)=b(1+y)^{-p}=(1+y)m(y)/p$. Hypothesis (W) holds with $\Lambda_{T_*}\le(1+T_*)^{\max\{2\theta-p-1,0\}}$, because
\[
 \frac{m(r+t)w(r+t)}{m(r)w(r)}=\Bigl(\frac{1+r+t}{1+r}\Bigr)^{2\theta-p-1}.
\]
The first moment of $m$ is finite because $p>1$. Here $\mathcal M_w=\mathcal M_\theta$ and $\mathscr X_w=\mathscr X_\theta$. By Propositions~\ref{prop:general-translation-norm} and~\ref{prop:abstract-weighted-transport}, there are $0<c'\le C'$ such that
\begin{equation}
 c'(1+t)^{-\alpha(p,\theta)}\le\rho_w(t)\le C'(1+t)^{-\alpha(p,\theta)},\quad\text{for all }t\ge0 .
 \label{eq:poly-rho}
\end{equation}
If $2\theta>1$, then $\alpha(p,\theta)>1$ and $\rho_w\in L^1(\R_+)$. Alternative (i) of Corollary~\ref{cor:general-two-sided} holds. If $2\theta\le1$, then $\alpha(p,\theta)=2\theta<p$ and
\[
 \kappa(t)=b(1+t)^{-p}\le b(1+t)^{-2\theta}\le\frac b{c'}\rho_w(t).
\]
Alternative (ii) holds. In both cases, Theorem~\ref{thm:general-transfer} and \eqref{eq:poly-rho} give $E(t)\le CN(1+t)^{-\alpha(p,\theta)}$ for every $U_0\in\mathscr X_\theta$, with $N=E(0)+\|\xi_0\|^2_{\mathcal M_\theta}$. This is \eqref{eq:sharp-abstract-transfer-upper}. Corollary~\ref{cor:general-two-sided} and \eqref{eq:poly-rho} give \eqref{eq:sharp-abstract-transfer-bilateral}.
\end{proof}

\begin{remark}[Operator norm and individual trajectories]
The lower estimate in \eqref{eq:sharp-abstract-transfer-bilateral} concerns the operator norm. Its proof uses pure-history initial states $(0,0,z_0)$, and the extremal history may depend on $t$. Individual trajectories are governed by their own inherited energy. Theorem~\ref{thm:prescribed-source-transfer} gives two-sided estimates for a fixed trajectory when its inherited energy satisfies matching upper and lower bounds.
\end{remark}

\begin{theorem}[Decay transfer from a prescribed history tail]
\label{thm:prescribed-source-transfer}
Assume the hypotheses of Theorem~\ref{thm:sharp-abstract-transfer},
except that the centered initial history is required only to belong to
$\mathcal M_m$. Let $U_0\in\mathscr X$, and suppose that for some
$\beta>0$ and $K\ge0$,
\begin{equation}
 H_{\xi_0}(t)=\int_0^\infty m(t+s)\|\xi_0(s)\|_Y^2\,ds
 \le K(1+t)^{-\beta},\quad\text{for all }t\ge0.
 \label{eq:prescribed-source-assumption}
\end{equation}
Then the corresponding energy satisfies
\begin{equation}
 E(t)\le C_{p,\beta}\bigl(E(0)+K\bigr)
 (1+t)^{-\min\{\beta,p+1\}},\quad\text{for all }t\ge0.
 \label{eq:prescribed-source-upper}
\end{equation}
The constant depends on the coefficients of the system, $p$ and
$\beta$, but not on $U_0$ or $K$. Furthermore,
\begin{equation}
 E(t)\ge c_*H_{\xi_0}(t),\quad\text{for all }t\ge0.
 \label{eq:prescribed-source-lower}
\end{equation}
If also $H_{\xi_0}(t)\ge k(1+t)^{-\beta}$ with $k>0$ and
$0<\beta\le p+1$, the trajectory itself has two-sided decay of
order $(1+t)^{-\beta}$.
\end{theorem}
\begin{proof}
The proof of Theorem~\ref{thm:sharp-abstract-transfer} shows that $m$ satisfies (K1)--(K3) with $a_*=p+1$, $C_m=1$ and $K_m=1/p$, and that its first moment is finite. Put $\lambda_*=\max\{p+1,\beta\}$ and $\sigma(t)=(1+t)^{-\beta}$. The functions $m$ and $\sigma$ have decay of order at most $\lambda_*$ with constant $1$, since $(1+t)/(1+s)\le t/s$ for $1\le s\le t$. Theorem~\ref{thm:source-transfer} applies with this $\sigma$, with $\lambda_*$ in place of $\lambda$. Estimate \eqref{eq:source-first} gives
\[
 E(t)\le CN\bigl(b(1+t)^{-p}+(1+t)^{-\beta}\bigr)\le C(b+1)N(1+t)^{-\min\{p,\beta\}} .
\]
If $\beta\le1$, then $\min\{p,\beta\}=\beta=\min\{\beta,p+1\}$, which is \eqref{eq:prescribed-source-upper}. If $\beta>1$, then $\sigma\in L^1(\R_+)$, and \eqref{eq:source-second} gives
\[
 E(t)\le C'N\bigl(pb(1+t)^{-p-1}+(1+t)^{-\beta}\bigr)\le C'(pb+1)N(1+t)^{-\min\{\beta,p+1\}} .
\]
This is \eqref{eq:prescribed-source-upper}. The lower estimate \eqref{eq:prescribed-source-lower} is part of Theorem~\ref{thm:source-transfer}, and combining both bounds proves the last assertion.
\end{proof}

\begin{remark}[Zero initial history and zero centered history]
\label{rem:zero-history-slow}
These conditions are different. If $z_0=0$ and
$\ell(u_0)\ne0$, then $\xi_0=-\ell(u_0)$ is constant and
\[
 H_{\xi_0}(t)=b\|\ell(u_0)\|_Y^2(1+t)^{-p}.
\]
Theorem~\ref{thm:prescribed-source-transfer} and its lower
history bound give matching energy estimates of order $(1+t)^{-p}$
for this trajectory. By contrast, $\xi_0=0$ means
$z_0(s)=\ell(u_0)$ almost everywhere and $H_{\xi_0}=0$.
The prescribed-source theorem gives an upper bound of order
$(1+t)^{-(p+1)}$. A matching lower bound for this trajectory requires
additional information on the observed motion.

For a continuous physical prehistory,
$z_0(s)=\ell(u_0)-\ell(u(-s))$.
Consequently $z_0=0$ represents an observation held at its
initial value, whereas $\xi_0=0$ represents a past observation
vanishing almost everywhere. If $\ell(u_0)\ne0$, the latter case
allows an initial discontinuity of the prescribed past at $t=0$
and is naturally interpreted as an energy-space initial state.
The compatible zero-prehistory regime of
\cite{AnasseBomissoKpata2026} also imposes $\ell(u_0)=0$.
\end{remark}

\begin{remark}[Weighted spaces and prescribed sources]
\label{rem:source-vs-weight}
Theorem~\ref{thm:sharp-abstract-transfer} follows as a special
case of the upper estimate in
Theorem~\ref{thm:prescribed-source-transfer}: for
$\xi_0\in\mathcal M_\theta$,
Proposition~\ref{prop:abstract-weighted-transport} gives
\eqref{eq:prescribed-source-assumption} with
$\beta=\min\{2\theta,p+1\}$ and
$K\le C\|\xi_0\|_{\mathcal M_\theta}^2$.
The prescribed-source condition also covers histories with
$\|\xi_0(s)\|_Y^2$ of order $(1+s)^r$, $0\le r<p$.
Their inherited source has order $(1+t)^{-(p-r)}$, while membership
in $\mathcal M_\theta$ holds for $2\theta<p-r$.
\end{remark}

\begin{corollary}[Strong stability from weighted decay]
\label{cor:strong-stability}
Under the hypotheses of Theorem~\ref{thm:sharp-abstract-transfer},
$T(t)U_0\to0$ in $\mathscr X$ for every $U_0\in\mathscr X$.
\end{corollary}
\begin{proof}
The union of the centered-history spaces $\mathscr X_\theta$
($\theta>0$) is dense in $\mathscr X$. Indeed, for fixed $u$,
$z-\ell(u)$ can be approximated in $\mathcal M_m$ by smooth,
compactly supported histories. Theorem~\ref{thm:sharp-abstract-transfer}
gives decay on this dense set. Since $T(t)$ is contractive on
$\mathscr X$, density extends the limit to every $U_0\in\mathscr X$.
\end{proof}

\begin{remark}[Saturation threshold and comparison of decay scales]
The exact saturation threshold in \eqref{eq:sharp-abstract-transfer-bilateral} is $\theta_c=(p+1)/2$. The value
$\theta=1/2$ is the threshold used in the bootstrap to obtain
$\varphi\in L^2(0,\infty)$.

The weighted-history scale measures decay through localization
of the past. It differs from the regularity-based scale described
by resolvent estimates such as those in \cite{BorichevTomilov}.
\end{remark}

\section{Applications}
\label{sec:other-applications}

We verify the generation and multiplier hypotheses for four mechanical realizations. Each uses the memory system \eqref{eq:variational-system}. The verifications of Assumption~\ref{ass:decay-transfer} below use the kernel only through its positivity, monotonicity, integrability and the value $m(0)$. Hence Corollary~\ref{cor:general-two-sided} applies to each realization for every pair $(m,w)$ satisfying its hypotheses. We state the explicit exponents for polynomial kernels through Theorem~\ref{thm:sharp-abstract-transfer}. The gradient-memory wave is the principal PDE application. Section~\ref{subsec:log-kernels} treats logarithmic kernels and weights for this model.

\subsection{Oscillators, strings, plates and gradient-memory waves}
\label{subsec:lower-order-examples}

\begin{example}[Damped oscillator with scalar memory]
\label{ex:oscillator}
Let $V=H=\C$, fix $\omega_0,\gamma>0$ and $\zeta\ge0$, and set
\[
 a(u,\psi)=\omega_0^2u\overline\psi,\qquad
 d(v,\psi)=\zeta v\overline\psi\quad\text{and}\quad\ell(u)=u.
\]
Equation~\eqref{eq:variational-system} becomes the system given by
\begin{equation*}
 \begin{cases}
 \dfrac{\partial^2u}{\partial t^2}+\omega_0^2u+\zeta\dfrac{\partial u}{\partial t}+\gamma\displaystyle\int_0^\infty m(s)z(t,s)\,ds=0,\\[2mm]
 \dfrac{\partial z}{\partial t}+\dfrac{\partial z}{\partial s}=\dfrac{\partial u}{\partial t},\\[2mm]
 z(t,0)=0.
 \end{cases}
\end{equation*}
Theorem~\ref{thm:variational-generation} gives generation. For $\zeta>0$, Proposition~\ref{prop:abstract-multiplier} verifies (H2). For $\zeta=0$, the coercive observation $\ell=\mathrm{id}$ permits Proposition~\ref{prop:full-memory-multiplier}. Thus Corollary~\ref{cor:general-two-sided} applies for all $\zeta\ge0$. For polynomial kernels, the weighted-history rate $\min\{2\theta,p+1\}$ and strong stability hold for all $\zeta\ge0$. The simulation below uses $\zeta=1$.
\end{example}

\begin{example}[One-dimensional wave equation with endpoint memory]
\label{ex:string}
Let $H=L^2(0,L)$ and $V=\{u\in H^1(0,L):u(0)=0\}$.
Fix $c,\eta,\gamma>0$ and $\zeta\ge0$, and put
\[
 a(u,\psi)=c^2\int_0^L u_x\overline{\psi_x}\,dx
 \quad\text{and}\quad \ell(u)=u(L),
\]
\[
 d(v,\psi)=\zeta\int_0^L v\overline\psi\,dx
 +\eta v(L)\overline{\psi(L)}.
\]
The weak equation describes
\begin{equation*}
 \begin{cases}
 \dfrac{\partial^2u}{\partial t^2}-c^2\dfrac{\partial^2u}{\partial x^2}+\zeta\dfrac{\partial u}{\partial t}=0,&0<x<L,\\
 u(t,0)=0,\\
 c^2\dfrac{\partial u}{\partial x}(t,L)+\eta\dfrac{\partial u}{\partial t}(t,L)+\gamma Jz(t)=0,\\
 \dfrac{\partial z}{\partial t}+\dfrac{\partial z}{\partial s}=\dfrac{\partial u}{\partial t}(t,L),\quad z(t,0)=0.
 \end{cases}
\end{equation*}
The form $a$ is coercive by the one-sided Poincar\'e inequality.
The trace $u\mapsto u(L)$ is continuous on $V$.
Also $\Real d(v,v)=\zeta\|v\|_H^2+\eta|v(L)|^2$ and
$|d(v,u)|\le C\bigl(\Real d(v,v)\bigr)^{1/2}a(u,u)^{1/2}$.
Generation therefore holds for every $\zeta\ge0$. If $\zeta>0$,
Proposition~\ref{prop:abstract-multiplier} and trace continuity give
(H1)--(H3), and Theorem~\ref{thm:sharp-abstract-transfer} applies.
Boundary-memory wave stabilization has been studied in related models,
for example in \cite{Mustafa2015,DuHao2024}.

We now treat the case $\zeta=0$, in which the damping acts only at
$x=L$ and $D_0(t)=\eta|u_t(t,L)|^2$. For a smooth real or complex solution,
set
\begin{equation*}
 \Psi_{\rm str}(t)=\Real\int_0^L x\,u_x(t,x)\overline{u_t(t,x)}\,dx.
\end{equation*}
Write $E_0=\frac12(\|u_t\|_{L^2}^2+c^2\|u_x\|_{L^2}^2)$.
We first identify the regularity of classical trajectories. Let $U=(u,v,z)\in D(\mathscr A_m)$. Then $v\in V\subset H^1(0,L)$, and there is $h\in L^2(0,L)$ such that
\begin{equation*}
 c^2\int_0^L u_x\overline{\psi_x}\,dx+\zeta\int_0^L v\overline\psi\,dx
 +\bigl(\eta v(L)+\gamma Jz\bigr)\overline{\psi(L)}
 =\int_0^L h\overline\psi\,dx,\quad\text{for all }\psi\in V.
\end{equation*}
Taking $\psi\in C_c^\infty(0,L)$ gives $-c^2u_{xx}=h-\zeta v$ in the sense of distributions. Hence $u_{xx}\in L^2(0,L)$ and $u\in H^2(0,L)$. Integration by parts with $\psi\in V$ and $\psi(L)=1$ then gives the boundary relation
\begin{equation*}
 c^2u_x(L)+\eta v(L)+\gamma Jz=0.
\end{equation*}
If $U_0\in D(\mathscr A_m)$, then $U(t)\in D(\mathscr A_m)$ and $U\in C^1([0,\infty);\mathscr X)$. Thus $u(t)\in H^2(0,L)$ and $u_t(t)\in H^1(0,L)$ for every $t\ge0$. Moreover $u\in C^1([0,\infty);H^1(0,L))$ and $u_t\in C^1([0,\infty);L^2(0,L))$, so $\Psi_{\rm str}$ is continuously differentiable. For these trajectories, integration by parts in $x$ gives the identity
\begin{equation*}
 \Psi_{\rm str}'(t)
 =-E_0(t)+\frac L2|u_t(t,L)|^2
       +\frac{Lc^2}{2}|u_x(t,L)|^2.
\end{equation*}
The memory boundary relation and weighted Cauchy--Schwarz imply
\[
 c^2|u_x(L)|\le\eta|u_t(L)|+\gamma|Jz|
 \quad\text{and}\quad |Jz|^2\le\kappa_0M,
\]
where $\kappa_0=\int_0^\infty m(s)\,ds$. Therefore
\begin{equation*}
 \Psi_{\rm str}'(t)\le-E_0(t)+C_\eta D_0(t)+C_{\eta,\gamma}M(t)
 \quad\text{and}\quad |\Psi_{\rm str}(t)|\le C E_0(t).
\end{equation*}
These estimates hold for classical trajectories. This proves (H2) even with $\zeta=0$, while the boundary dissipation
and the string trace inequality give (H1) and (H3).
Consequently the sharp two-sided weighted-history rate applies to the
string with damping solely at the observed endpoint.
\end{example}

\begin{example}[Clamped plate with rank-one distributed memory]
\label{ex:plate}
Let $\Omega\subset\R^2$ be bounded and smooth,
$H=L^2(\Omega)$, and $V=H_0^2(\Omega)$.
Fix $\rho\in L^2(\Omega)$ with $\rho\ne0$. The forms and the observation are given by
\[
 \begin{aligned}
 a(u,\psi)&=\int_\Omega\Delta u\,\overline{\Delta\psi}\,dx,\qquad
 d(v,\psi)=\zeta\int_\Omega v\overline\psi\,dx,\\
 \text{and}\quad\ell(u)&=\int_\Omega\rho(x)u(x)\,dx,
 \end{aligned}
\]
where $\rho$ is taken real valued for simplicity and $\zeta>0$.
The associated system is
\begin{equation*}
 \begin{cases}
 \dfrac{\partial^2u}{\partial t^2}+\Delta^2u+\zeta\dfrac{\partial u}{\partial t}+\gamma\rho(x)Jz=0,
 &x\in\Omega,\\
 u=\partial_\nu u=0,&x\in\partial\Omega,\\
 \dfrac{\partial z}{\partial t}+\dfrac{\partial z}{\partial s}=\displaystyle\int_\Omega\rho(x)\dfrac{\partial u}{\partial t}(x,t)\,dx,\\
 z(t,0)=0.
 \end{cases}
\end{equation*}
Standard elliptic estimates make $a$ coercive on $H_0^2(\Omega)$, while $\ell$ is bounded on $H$. With $\zeta>0$, Theorem~\ref{thm:variational-generation} and Proposition~\ref{prop:abstract-multiplier} give the sharp weighted decay rate. Viscoelastic plate models with infinite memory are studied in \cite{AnayaMessaoudiMustapha2021}. If $\zeta=0$ and an eigenmode $e_j$ satisfies $\ell(e_j)=0$, then $u(t)=\cos(\sqrt{\lambda_j}t)e_j$, $z=0$ is undamped. Moreover, $\ell(e_j)\to0$ for every orthonormal eigenbasis, so the memory observation does not control the high modes uniformly. Our verification of (H2) therefore assumes $\zeta>0$.
\end{example}

\begin{example}[Wave equation with distributed gradient memory]
\label{ex:gradient-wave}
Let $\Omega\subset\R^d$ be a bounded smooth domain, let
$H=L^2(\Omega)$, $V=H_0^1(\Omega)$, and put
\[
 Y=\operatorname{Ran}(\nabla|_{H_0^1(\Omega)})\subset L^2(\Omega)^d
 \quad\text{and}\quad \ell(u)=\nabla u.
\]
The space $Y$ is closed, because the Poincar\'e inequality makes
$\nabla:H_0^1(\Omega)\to L^2(\Omega)^d$ bounded below when
$H_0^1(\Omega)$ carries its usual gradient norm. Choose $k>0$, $\zeta\ge0$ and $\gamma>0$, and set
\[
 a(u,\psi)=k\int_\Omega\nabla u\cdot\overline{\nabla\psi}\,dx
 \quad\text{and}\quad
 d(v,\psi)=\zeta\int_\Omega v\overline\psi\,dx.
\]
The abstract formulation becomes
\begin{equation*}
 \begin{cases}
 \dfrac{\partial^2u}{\partial t^2}-k\Delta u+\zeta\dfrac{\partial u}{\partial t}-\gamma\operatorname{div}Jz=0,
 &x\in\Omega,\ t>0,\\
 u(t)|_{\partial\Omega}=0,\\
 \dfrac{\partial z}{\partial t}(t,s)+\dfrac{\partial z}{\partial s}(t,s)=\nabla\dfrac{\partial u}{\partial t}(t)\quad\text{and}\quad z(t,0)=0,\\
 Jz(t,x)=\displaystyle\int_0^\infty m(s)z(t,s,x)\,ds.
 \end{cases}
\end{equation*}
The divergence is understood in the distributional sense, so
$\ell^*y=-\operatorname{div}y\in H^{-1}(\Omega)$.
If $\zeta>0$, then
$\Real d(v,v)=\zeta\|v\|_{L^2}^2$ and
\[
 |d(v,u)|\le\zeta\|v\|_{L^2}\|u\|_{L^2}
 \le C\bigl(\Real d(v,v)\bigr)^{1/2}a(u,u)^{1/2}.
\]
Thus Proposition~\ref{prop:abstract-multiplier} verifies (H2).
When $\zeta=0$, no instantaneous damping is present. Nevertheless,
by the Poincar\'e inequality, the observation
$\ell=\nabla: H_0^1(\Omega)\to Y$ is bounded below. With the gradient
norm on $V$ it is an isometric isomorphism onto $Y$. Consequently,
Proposition~\ref{prop:full-memory-multiplier} proves
(H2) even for $\zeta=0$. The lifting is $\eta=\ell^{-1}z\in V$.
In both cases, continuity of $\nabla$ gives (H3), generation follows
from Theorem~\ref{thm:variational-generation}, and the invariance of the weighted spaces
from Proposition~\ref{prop:general-invariance}.
Therefore Corollary~\ref{cor:general-two-sided} applies for all $\zeta\ge0$ and every pair $(m,w)$ satisfying its hypotheses. For polynomial kernels and weights, Theorem~\ref{thm:sharp-abstract-transfer} gives the optimal exponent $\min\{2\theta,p+1\}$.

For $\zeta=0$, put $\kappa_0=\int_0^\infty m(s)\,ds$ and $g(s)=\gamma m(s)$. The characteristic formula gives
\[
 Jz(t)=\kappa_0\nabla u(t)-\int_0^\infty m(s)\nabla u(t-s)\,ds,
\]
where the negative-time displacement is prescribed by the initial history. If $k=1-\gamma\kappa_0>0$, the displacement equation is exactly
\begin{equation*}
 \begin{cases}
 \dfrac{\partial^2u}{\partial t^2}-\Delta u+\displaystyle\int_0^\infty g(s)\Delta u(t-s)\,ds=0, & x\in\Omega,\ t>0,\\[2mm]
 u=0, & x\in\partial\Omega,\ t>0,
 \end{cases}
 \quad\text{with}\quad \int_0^\infty g(s)\,ds=\gamma\kappa_0<1 .
\end{equation*}
On a Dirichlet eigenfunction $-\Delta e_j=\lambda_je_j$, the displacement amplitude $x_j$ satisfies
\[
 x_j''(t)+\lambda_jx_j(t)-\lambda_j\int_0^\infty g(s)x_j(t-s)\,ds=0 .
\]
Indeed, the augmented equation has elastic coefficient $1-\gamma\kappa_0$, and its history force adds $\gamma\kappa_0\lambda_jx_j(t)$ and subtracts the same convolution. The energy is given by
\[
 \mathcal E(t)=\frac12\Bigl[\|u_t(t)\|_{L^2}^2+(1-\gamma\kappa_0)\|\nabla u(t)\|_{L^2}^2+\gamma\int_0^\infty m(s)\|\nabla u(t)-\nabla u(t-s)\|_{L^2}^2\,ds\Bigr].
\]
If $z_0$ arises from a prescribed displacement history $u(-s)$, then $\xi_0(s)=-\nabla u(-s)$ for $s>0$. The weighted norm of $\xi_0$ therefore measures the localization of the prescribed past. For $m(s)=pb(1+s)^{-p-1}$, so that $\kappa_0=b$, this equation is a polynomial-kernel subclass of the first infinite-memory wave model in \cite{Guesmia2020}, with $g(s)=\gamma pb(1+s)^{-p-1}$ and exponent $q=p+1>2$. Corollary~\ref{cor:log-waves} below gives the two-sided law for this model without instantaneous damping.
\end{example}

\begin{corollary}[Algebraically growing past histories]
\label{cor:algebraic-history-data}
Let $m(s)=pb(1+s)^{-p-1}$, with $p>1$ and $b>0$. Suppose the abstract system satisfies
Assumption~\ref{ass:decay-transfer}, and let $0\le r<p$.
If the centered initial history satisfies
\begin{equation*}
 \|\xi_0(s)\|_Y^2\le K_0(1+s)^r
 \quad\text{for almost every }s\ge0,
\end{equation*}
then
\begin{equation*}
 E(t)\le C_{p,r}\bigl(E(0)+K_0\bigr)
 (1+t)^{-(p-r)},\quad\text{for all }t\ge0.
\end{equation*}
If additionally, for some $k_0>0$ and $s_0\ge0$,
\begin{equation*}
 \|\xi_0(s)\|_Y^2\ge k_0(1+s)^r
 \quad\text{for almost every }s\ge s_0,
\end{equation*}
then
\begin{equation*}
 E(t)\ge c_{p,r}k_0(1+t)^{-(p-r)}
 \quad\text{for }t\ge\max\{1,s_0\}.
\end{equation*}
Thus the exponent $p-r$ is attained by every trajectory with
these two-sided history bounds. The estimate is a pure power law with no logarithmic
factor in the upper estimate, including when $p-r=1$.
\end{corollary}
\begin{proof}
The assumption $r<p$ ensures that $\xi_0\in\mathcal M_m$.
Moreover, because $1+s\le1+t+s$ and $m(s)=pb(1+s)^{-p-1}$,
\[
 \begin{aligned}
 H_{\xi_0}(t)
 &\le pb K_0\int_0^\infty
   (1+t+s)^{-p-1}(1+s)^r\,ds\\
 &\le pb K_0\int_0^\infty
   (1+t+s)^{-p-1+r}\,ds\\
 &=\frac{pb}{p-r}K_0(1+t)^{-(p-r)}.
 \end{aligned}
\]
Theorem~\ref{thm:prescribed-source-transfer} applies with
$\beta=p-r$ and proves the upper bound. Under the additional
lower growth assumption, for $t\ge\max\{1,s_0\}$ we have
\[
 \begin{aligned}
 H_{\xi_0}(t)
 &\ge pb k_0\int_t^{2t}
      (1+t+s)^{-p-1}(1+s)^r\,ds\\
 &\ge pb k_0\,t(1+3t)^{-p-1}(1+t)^r\\
 &\ge c_{p,r}k_0(1+t)^{-(p-r)}.
 \end{aligned}
\]
The lower estimate follows from
\eqref{eq:prescribed-source-lower}.
\end{proof}

\begin{remark}[Comparison with Guesmia's polynomial-kernel example]
\label{rem:guesmia-prior-comparison}
Consider the first viscoelastic wave model of
Guesmia~\cite{Guesmia2020}, Example~1, with
$g(s)=a(1+s)^{-q}$ and the two-sided history condition given by
\[
 m_1(1+s)^r\le1+\|\nabla u(-s)\|_{L^2}^2
 \le m_0(1+s)^r\quad\text{with}\quad 0\le r<q-1.
\]
In our polynomial subclass,
$q=p+1>2$, $a=\gamma pb$ and $a<q-1$ is equivalent to
$\gamma b<1$. The centered initial history is
$\xi_0(s)=-\nabla u(-s)$. In formula (4.28) of
\cite{Guesmia2020}, the energy is bounded by a constant times
$(1+t)^{-(q-r-1)}$ if $q-r\ne2$, and by
$(1+t)^{-1}(1+\log(1+t))$ if $q-r=2$.
Corollary~\ref{cor:algebraic-history-data} gives, for the
polynomial subclass of the first model,
\[
 E(t)\le C(1+t)^{-(q-r-1)}
\]
\emph{also when $q-r=2$}, with no logarithmic factor. At the borderline $q-r=2$, one has $r=p-1>0$.
If $r>0$, Guesmia's lower growth hypothesis implies
$\|\xi_0(s)\|_Y^2\ge c(1+s)^r$ for all sufficiently large
$s$, so our corollary supplies the matching individual lower
bound. For $r=0$, the same lower estimate holds under the explicit
condition $\|\xi_0(s)\|_Y^2\ge c_0>0$ for all sufficiently large $s$.
\end{remark}

\subsection{Logarithmic kernels and logarithmic weights}
\label{subsec:log-kernels}

We apply Corollary~\ref{cor:general-two-sided} to a kernel and a weight that both carry a logarithmic factor. The two factors are independent. Let $p>1$, $\alpha\in\R$, $\eta\ge0$, $\theta>0$ and $c>0$. We write $\alpha_+=\max\{\alpha,0\}$, $\alpha_-=\max\{-\alpha,0\}$, $L(s)=\log(e+s)$ and $T=1+t$. We assume the following condition:
\begin{equation}
 \alpha_+<p .
 \label{eq:log-alpha}
\end{equation}
The kernel and the weight are given by
\begin{equation}
 m(s)=c(1+s)^{-p-1}L(s)^{\alpha}\quad\text{and}\quad w(s)=(1+s)^{2\theta}L(s)^{\eta},\quad\text{for all }s\ge0 .
 \label{eq:log-kernel}
\end{equation}
Since $L(0)=1$, the translation norm \eqref{eq:rho-w} is given by
\begin{equation}
 \rho_w(t)=\sup_{r\ge0}Q(r,t)\quad\text{with}\quad Q(r,t)=\frac{(1+r)^{p+1-2\theta}L(r+t)^{\alpha}}{(1+r+t)^{p+1}L(r)^{\alpha+\eta}} .
 \label{eq:log-rho}
\end{equation}
The supremum replaces the essential supremum because $Q$ is continuous in $r$. We use the elementary inequality
\begin{equation}
 \log y\le\frac{y^{\beta}}{e\beta},\quad\text{for all }y>0\text{ and }\beta>0,
 \label{eq:log-pow}
\end{equation}
and the increment bound
\begin{equation}
 1\le\frac{L(r+t)}{L(r)}\le1+\log\Bigl(1+\frac t{1+r}\Bigr),\quad\text{for all }r,t\ge0 .
 \label{eq:log-increment}
\end{equation}
The latter follows from $L(r+t)-L(r)=\log(1+t/(e+r))$ and $L(r)\ge1$.

\begin{proposition}[Structural hypotheses]
\label{prop:log-hypotheses}
Under \eqref{eq:log-alpha}, the pair \eqref{eq:log-kernel} satisfies (K1)--(K3) and (W) with
\[
 a_*=p+1-\alpha_+,\quad\lambda=p+1+\alpha_-,\quad C_m=1,\quad K_m=\frac1{p-\alpha_+}\quad\text{and}\quad\Lambda_{T_*}\le(1+T_*)^{|p+1-2\theta|+|\alpha+\eta|}.
\]
Moreover $\int_0^\infty s\,m(s)\,ds<\infty$, and $m$ satisfies the hypotheses of Theorem~\ref{thm:variational-generation} and Proposition~\ref{prop:full-memory-multiplier}.
\end{proposition}
\begin{proof}
\emph{(K1).} The logarithmic derivative is given by
\[
 -\frac{m'(s)}{m(s)}=\frac{p+1}{1+s}-\frac{\alpha}{(e+s)L(s)} .
\]
Since $1+s\le(e+s)L(s)$, the last term is at most $\alpha_+/(1+s)$. Hence $-m'(s)\ge a_*m(s)/(1+s)$, with $a_*>1$ by \eqref{eq:log-alpha}. In particular $m$ is smooth, positive, decreasing, integrable and $m(0)=c$.

\emph{(K2).} Let $1\le s\le t$. Then $(1+t)/(1+s)\le t/s$. If $\alpha\ge0$, then $m(s)/m(t)\le(t/s)^{p+1}$. If $\alpha<0$, \eqref{eq:log-increment} gives $L(t)/L(s)\le1+\log(t/s)\le t/s$, hence $m(s)/m(t)\le(t/s)^{p+1+\alpha_-}$.

\emph{(K3).} Integration of (K1) gives $m(x)\le m(y)\bigl((1+y)/(1+x)\bigr)^{a_*}$ for $x\ge y\ge0$. Since $a_*>1$,
\[
 \kappa(y)\le m(y)(1+y)^{a_*}\int_y^\infty(1+x)^{-a_*}\,dx=\frac{(1+y)m(y)}{a_*-1} .
\]

\emph{(W).} The weight satisfies $w\ge1$, since $\eta\ge0$, and $mw\in L^1_{\rm loc}$. The product is given by $mw=c(1+s)^{2\theta-p-1}L(s)^{\alpha+\eta}$. For $0\le t\le T_*$, both $(1+r+t)/(1+r)$ and $L(r+t)/L(r)$ lie in $[1,1+T_*]$. This gives the bound on $\Lambda_{T_*}$.

\emph{First moment.} If $\alpha>0$, \eqref{eq:log-pow} with $\beta=(p-1)/(2\alpha)$ gives $L(s)^{\alpha}\le C(1+s)^{(p-1)/2}$. Hence $s\,m(s)\le C(1+s)^{-(p+1)/2}$, which is integrable because $p>1$. If $\alpha\le0$, then $s\,m(s)\le c(1+s)^{-p}$.
\end{proof}

\begin{theorem}[Explicit translation rates]
\label{thm:log-rates}
Assume \eqref{eq:log-alpha} and $\eta\ge0$. Put $\nu=\alpha_++\eta$. The following estimates hold for every $t\ge0$.
\begin{enumerate}[label=\textup{(\roman*)},leftmargin=*]
\item If $2\theta<p+1$, then, with
\[
 c_1=2^{-p-1}(1+\log2)^{-\alpha_-}\quad\text{and}\quad C_1=2^{p+1}\Bigl(1+\frac{2\nu}{e(p+1-2\theta)}\Bigr)^{\nu},
\]
we have
\[
 c_1T^{-2\theta}L(t)^{-\eta}\le\rho_w(t)\le C_1T^{-2\theta}L(t)^{-\eta}.
\]
\item If $2\theta>p+1$ and $\alpha\ge0$, then $\rho_w(t)=T^{-p-1}L(t)^{\alpha}$. If $2\theta>p+1$ and $\alpha<0$, put $\mu=(\alpha_--\eta)_+$. Then
\[
 T^{-p-1}L(t)^{\alpha}\le\rho_w(t)\le C_2T^{-p-1}L(t)^{\alpha},
\]
with $C_2=1$ if $\mu=0$, and $C_2=(e^{\beta-1}/\beta)^{\mu}$ with $\beta=(2\theta-p-1)/\mu$ if $\mu>0$.
\item If $2\theta=p+1$ and $\alpha\ge0$, then $\rho_w(t)=T^{-p-1}L(t)^{\alpha}$. If $2\theta=p+1$ and $\alpha<0$, then
\[
 2^{-p-1}(1+\log2)^{\alpha}\,T^{-p-1}L(t)^{-\min\{\alpha_-,\eta\}}\le\rho_w(t)\le T^{-p-1}L(t)^{-\min\{\alpha_-,\eta\}} .
\]
\end{enumerate}
\end{theorem}
\begin{proof}
We bound the quotient $Q$ of \eqref{eq:log-rho}, and we write $u=1+r\ge1$. Then
\begin{equation}
 1+r+t\ge\frac{u+T}2\quad\text{and}\quad1+r+t\ge T .
 \label{eq:log-half}
\end{equation}

\emph{Lower bounds.} By \eqref{eq:log-increment} with $r=t$, we have $1\le L(2t)/L(t)\le1+\log2$. Since $1+2t\le2T$,
\begin{equation}
 \rho_w(t)\ge Q(t,t)=\frac{T^{p+1-2\theta}}{(1+2t)^{p+1}}\Bigl(\frac{L(2t)}{L(t)}\Bigr)^{\alpha}L(t)^{-\eta}\ge c_1T^{-2\theta}L(t)^{-\eta}.
 \label{eq:log-rt}
\end{equation}
Since $L(0)=1$,
\begin{equation}
 \rho_w(t)\ge Q(0,t)=T^{-p-1}L(t)^{\alpha}.
 \label{eq:log-r0}
\end{equation}

\emph{Upper bound in (i).} Since $L(t)\le L(r+t)$ and $L(r+t)/L(r)\ge1$,
\[
 Q(r,t)\,T^{2\theta}L(t)^{\eta}\le T^{2\theta}\frac{u^{p+1-2\theta}}{(1+r+t)^{p+1}}\Bigl(\frac{L(r+t)}{L(r)}\Bigr)^{\nu}.
\]
Put $x=u/T$. Since $t/(1+r)\le1/x$, \eqref{eq:log-increment} and \eqref{eq:log-half} give
\[
 Q(r,t)\,T^{2\theta}L(t)^{\eta}\le2^{p+1}\frac{x^{p+1-2\theta}}{(1+x)^{p+1}}\Bigl(1+\log\bigl(1+\tfrac1x\bigr)\Bigr)^{\nu}.
\]
If $\nu=0$, the right-hand side is at most $2^{p+1}(1+x)^{-2\theta}\le2^{p+1}$. If $\nu>0$, set $\beta=(p+1-2\theta)/(2\nu)$. By \eqref{eq:log-pow},
\[
 1+\log\bigl(1+\tfrac1x\bigr)\le\Bigl(1+\frac1{e\beta}\Bigr)\Bigl(1+\frac1x\Bigr)^{\beta}.
\]
Since $\nu\beta=(p+1-2\theta)/2$, we have
\[
 x^{p+1-2\theta}\Bigl(1+\frac1x\Bigr)^{(p+1-2\theta)/2}=x^{(p+1-2\theta)/2}(1+x)^{(p+1-2\theta)/2}\le(1+x)^{p+1-2\theta}.
\]
The right-hand side is therefore at most $2^{p+1}(1+1/(e\beta))^{\nu}(1+x)^{-2\theta}\le C_1$. The lower bound is \eqref{eq:log-rt}.

\emph{Cases (ii) and (iii) with $\alpha\ge0$.} The function defined by $\phi(y)=y^{-p-1}\log(e+y-1)^{\alpha}$, $y\ge1$, satisfies
\[
 \frac{\phi'(y)}{\phi(y)}=-\frac{p+1}y+\frac{\alpha}{(e+y-1)\log(e+y-1)}\le\frac{\alpha-p-1}y<0 .
\]
Since $1+r+t\ge T$, we get $L(r+t)^{\alpha}(1+r+t)^{-p-1}=\phi(1+r+t)\le\phi(T)=T^{-p-1}L(t)^{\alpha}$. Moreover $u^{p+1-2\theta}\le1$ and $L(r)^{-\alpha-\eta}\le1$. Hence $Q(r,t)\le T^{-p-1}L(t)^{\alpha}$. Together with \eqref{eq:log-r0}, this gives the identity.

\emph{Case (ii) with $\alpha<0$.} Since $T\le1+r+t$ and $L(t)\le L(r+t)$,
\[
 Q(r,t)\,T^{p+1}L(t)^{\alpha_-}\le u^{p+1-2\theta}\Bigl(\frac{L(t)}{L(r+t)}\Bigr)^{\alpha_-}L(r)^{\alpha_--\eta}\le u^{p+1-2\theta}L(r)^{\mu}.
\]
If $\mu=0$, this is at most $1$. If $\mu>0$, \eqref{eq:log-pow} and $e+r\le eu$ give $L(r)\le e^{\beta-1}u^{\beta}/\beta$. Since $\beta\mu=2\theta-p-1$, we get $u^{p+1-2\theta}L(r)^{\mu}\le C_2$. The lower bound is \eqref{eq:log-r0}.

\emph{Case (iii) with $\alpha<0$.} Here
\[
 Q(r,t)=\frac{L(r)^{\alpha_--\eta}}{(1+r+t)^{p+1}L(r+t)^{\alpha_-}} .
\]
If $\eta\ge\alpha_-$, then $L(r)^{\alpha_--\eta}\le1$ and $L(r+t)\ge L(t)$. Hence $Q(r,t)\le T^{-p-1}L(t)^{-\alpha_-}$. If $\eta<\alpha_-$, then $L(r)^{\alpha_--\eta}\le L(r+t)^{\alpha_--\eta}$. Hence $Q(r,t)\le T^{-p-1}L(r+t)^{-\eta}\le T^{-p-1}L(t)^{-\eta}$. In both cases $Q(r,t)\le T^{-p-1}L(t)^{-\min\{\alpha_-,\eta\}}$. For the lower bound, \eqref{eq:log-r0} gives $T^{-p-1}L(t)^{-\alpha_-}$ and \eqref{eq:log-rt} gives $2^{-p-1}(1+\log2)^{\alpha}T^{-p-1}L(t)^{-\eta}$. The larger of the two is at least $2^{-p-1}(1+\log2)^{\alpha}T^{-p-1}L(t)^{-\min\{\alpha_-,\eta\}}$.
\end{proof}

The roles of the two logarithms are separated by the threshold $2\theta=p+1$. Below it, the rate is $T^{-2\theta}L(t)^{-\eta}$. The logarithm of the weight enters and the logarithm of the kernel does not. Above it, the rate is $m(t)/c$. The logarithm of the kernel enters and the logarithm of the weight does not. At the threshold with $\alpha\ge0$ the kernel decides, and with $\alpha<0$ the smaller logarithmic gain decides. These rates depend on the pair $(m,w)$ and cannot be read on the kernel alone.

\begin{corollary}[Viscoelastic waves with logarithmic kernels and weights]
\label{cor:log-waves}
Let $\Omega\subset\R^d$ be a bounded smooth domain. Let $m$ and $w$ be given by \eqref{eq:log-kernel} with $p>1$, $\eta\ge0$ and \eqref{eq:log-alpha}. Put $\kappa_0=\int_0^\infty m(s)\,ds$, and let $\gamma>0$ satisfy $\gamma\kappa_0<1$. Consider the system of Example~\ref{ex:gradient-wave} with $\zeta=0$ and $k=1-\gamma\kappa_0$. Its displacement equation is
\[
 \begin{cases}
 \dfrac{\partial^2u}{\partial t^2}-\Delta u+\gamma\displaystyle\int_0^\infty m(s)\Delta u(t-s)\,ds=0, & x\in\Omega,\ t>0,\\[2mm]
 u=0, & x\in\partial\Omega,\ t>0.
 \end{cases}
\]
There are constants $0<c_W\le C_W$ such that
\[
 c_W\rho_w(t)\le\|T(t)\|^2_{\mathcal L(\mathscr X_w,\mathscr X)}\le C_W\rho_w(t),\quad\text{for all }t\ge0,
\]
where $\rho_w$ obeys Theorem~\ref{thm:log-rates}. In particular, for $\alpha=\eta=0$ and $c=pb$, there are $0<c_\theta\le C_\theta$ such that
\[
 c_\theta(1+t)^{-\min\{2\theta,p+1\}}\le\|T(t)\|^2_{\mathcal L(\mathscr X_\theta,\mathscr X)}\le C_\theta(1+t)^{-\min\{2\theta,p+1\}},\quad\text{for all }t\ge0 .
\]
This two-sided estimate holds without an instantaneous velocity damping term.
\end{corollary}
\begin{proof}
Example~\ref{ex:gradient-wave} and Proposition~\ref{prop:full-memory-multiplier} give Assumption~\ref{ass:decay-transfer} with $\zeta=0$, since $m$ is positive, decreasing, in $W^{1,1}(\R_+)$ and $m(0)<\infty$. Proposition~\ref{prop:log-hypotheses} gives (K1)--(K3), (W) and the finiteness of the first moment.

If $2\theta\ge p+1$, Theorem~\ref{thm:log-rates} gives $\rho_w(t)\le CT^{-p-1}L(t)^{\alpha_+}$, which is integrable. If $1<2\theta<p+1$, it gives $\rho_w(t)\le C_1T^{-2\theta}$, which is integrable. In these cases alternative (i) of Corollary~\ref{cor:general-two-sided} holds.

If $2\theta\le1$, (K3) gives $\kappa(t)\le cK_mT^{-p}L(t)^{\alpha}$. Theorem~\ref{thm:log-rates}(i) gives $\rho_w(t)\ge c_1T^{-2\theta}L(t)^{-\eta}$. By \eqref{eq:log-pow}, $L(t)^{\alpha_++\eta}\le CT^{(p-2\theta)/2}$. Hence
\[
 \frac{\kappa(t)}{\rho_w(t)}\le\frac{cK_m}{c_1}T^{-(p-2\theta)}L(t)^{\alpha_++\eta}\le C,
\]
because $p-2\theta\ge p-1>0$. Alternative (ii) holds. Corollary~\ref{cor:general-two-sided} concludes.
\end{proof}

\subsection{Numerical illustration: finite-duration initial history}\label{subsec:oscillator-numerics}

\begin{figure}[htbp]
\centering
\includegraphics[width=0.92\linewidth]{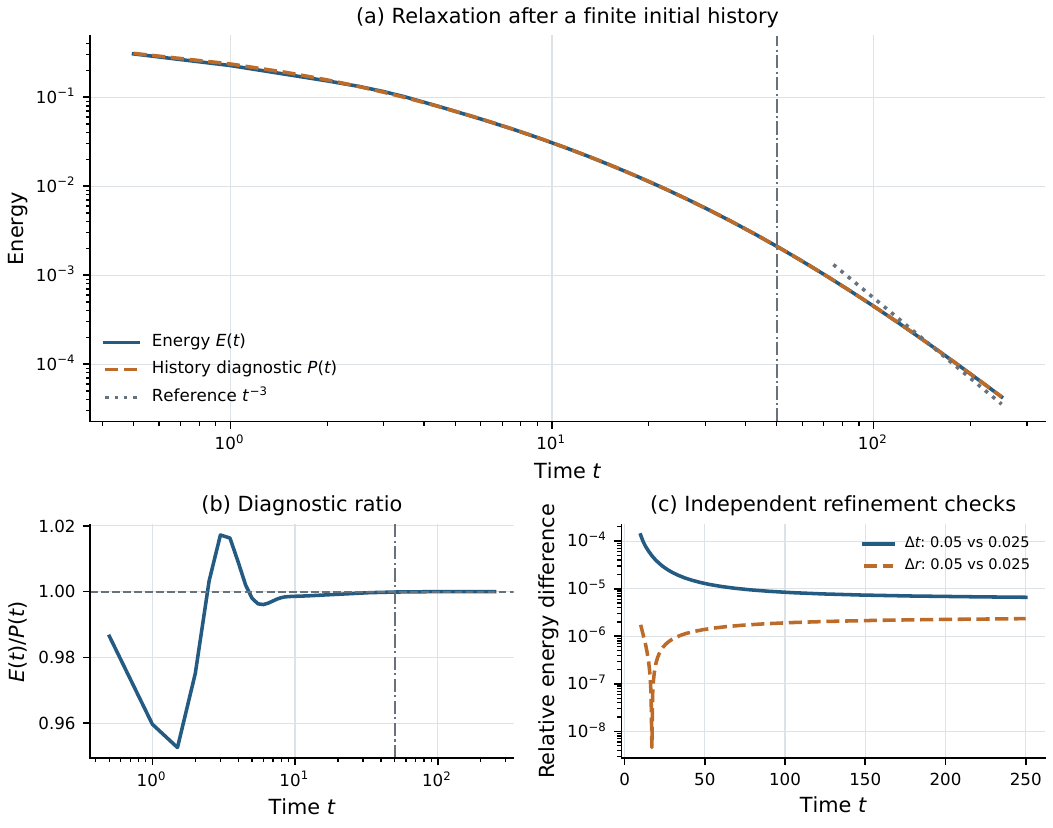}
\caption{Scalar oscillator with $\omega_0=\zeta=\gamma=b=1$, $p=2$, $T=50$, $\Delta t=\Delta r=0.05$. (a) Complete numerical energy $E(t)$, historical diagnostic $P(t)$, and a $t^{-3}$ reference guide. The vertical line marks $T$. (b) Ratio $E(t)/P(t)$. (c) Relative energy changes under separate halving of the time step and inherited-history quadrature step.}
\label{fig:oscillator-history-recovery}
\end{figure}

For Example~\ref{ex:oscillator}, choose $\omega_0=\zeta=\gamma=b=1$, $p=2$, $u(0)=u_t(0)=0$, and the finite-duration centered history defined by
\[
 \xi_0(r)=(1-e^{-r})(1+r)^{1/2}
 \left[\max\left\{0,\min\left(1,\frac{T-r}{10}\right)\right\}\right]^2
 \quad\text{with}\quad T=50.
\]
It vanishes at $r=0$ and for $r\ge T$ and belongs to every weighted history space. The characteristics reduce the coupled system to the scalar Volterra equation
\begin{equation*}
 u_{tt}+\zeta u_t+(\omega_0^2+\gamma b)u
 =\gamma\int_0^t m(s)u(t-s)\,ds-\gamma F_{\xi_0}(t),
\end{equation*}
with
\begin{equation*}
 F_{\xi_0}(t)=\int_0^T m(t+r)\xi_0(r)\,dr.
\end{equation*}
We compute the convolution by causal composite quadrature and advance the oscillator by a first-order semi-implicit time step. The energy is evaluated from the complete characteristic decomposition. In particular, although $\xi_0(r)=0$ for $r>T$, one has $z(t,t+r)=u(t)$ there. The inherited contribution is therefore
\[
 \int_0^T m(t+r)|u(t)+\xi_0(r)|^2\,dr
 +b(1+t+T)^{-p}|u(t)|^2.
\]
The first term is evaluated by quadrature and the second analytically. The comparison expression is defined as follows:
\[
 P(t)=\frac{\gamma}{2}\left[H_{\xi_0}(t)
 +m(t)\int_0^t|u(r)|^2\,dr\right].
\]
Figure~\ref{fig:oscillator-history-recovery} compares $E(t)$ with the historical diagnostic $P(t)$. For this nonzero, compactly supported history, the saturated rate is also sharp along the particular trajectory, not merely in operator norm. Indeed, by dominated convergence,
\[
 \lim_{t\to\infty}(1+t)^{p+1}H_{\xi_0}(t)
 =pb\int_0^T|\xi_0(r)|^2\,dr>0.
\]
The function $(1+t)^{p+1}H_{\xi_0}(t)$ is continuous and strictly positive for $t\ge0$. Its positive limit therefore gives a positive lower bound on the entire half-line. As $u(0)=u_t(0)=0$ and $\xi_0\in\mathcal M_\theta$ for every $\theta>0$, the upper bound in Theorem~\ref{thm:sharp-abstract-transfer}, with $\theta\ge(p+1)/2$, and the lower bound \eqref{eq:prescribed-source-lower} yield constants $c_{\xi_0},C_{\xi_0}>0$ such that
\[
 c_{\xi_0}(1+t)^{-(p+1)}\le E(t)
 \le C_{\xi_0}(1+t)^{-(p+1)},\quad\text{for all }t\ge0.
\]
The plot concerns the saturated rate. The non-saturated regime is tested separately by evaluating the exact translation norm in Proposition~\ref{prop:abstract-weighted-transport}, with histories localized near the maximizing age $s_t$.

When the time step is reduced from $0.05$ to $0.025$, the relative changes in the computed energy at $t=50$, $100$, and $200$ are $1.3\times10^{-5}$, $8.4\times10^{-6}$, and $6.9\times10^{-6}$, respectively. Holding the time step fixed at $0.05$ and independently halving the quadrature step for the inherited history from $0.05$ to $0.025$ gives relative energy changes of $1.4\times10^{-6}$, $1.9\times10^{-6}$, and $2.3\times10^{-6}$ at the same times. The source-history integrals are evaluated independently of the mechanical time-stepping grid in this check.

\paragraph{Non-saturated weighted-history diagnostic.}
The exact translation norm gives a numerical test of the feature not visible in the compact-history oscillator experiment. For $p=2$ and $\theta=1/2$, the maximizer of \eqref{eq:weighted-translation-exact} is $1+s_t=2t$, for $t\ge1$. Consequently
\[
 \|R(t)\|^2_{\mathcal L(\mathcal M_{1/2},\mathcal M_m)}
 =\frac{(2t)^2}{(3t)^3}=\frac{4}{27t},\quad\text{for all }t\ge1.
\]
For $\theta=1$ and $t\ge2$, the maximizer is $1+s_t=t/2$, yielding
\[
 \|R(t)\|^2_{\mathcal L(\mathcal M_1,\mathcal M_m)}
 =\frac{4}{27t^2}.
\]
For $\theta\ge3/2$, the maximum occurs at $s=0$, so the exact squared norm is $(1+t)^{-3}$. These three cases display the slopes $-1$, $-2$ and $-3$ and their saturation.

To evaluate the non-saturated norms without relying on the optimizer formula, take the centered history
\[
 \xi_{t,\theta}(s)=\begin{cases}
 \bigl(1-((s-s_t)/w_t)^2\bigr)^2,& |s-s_t|<w_t,\\
 0,&\text{otherwise},\end{cases}
 \quad\text{with}\quad w_t=0.025t,
\]
where $s_t$ is the maximizing age given above. For $t\in\{20,40,80,160\}$ this history is supported in $(0,\infty)$. We evaluate the normalized inherited energy
\[
 Q_\theta(t)=\frac{H_{\xi_{t,\theta}}(t)}{\|\xi_{t,\theta}\|_{\mathcal M_\theta}^2}
\]
by a 96-point Gauss--Legendre rule. The results are
\[
\begin{array}{c|cc|cc}
 &\multicolumn{2}{c|}{\theta=1/2}&\multicolumn{2}{c}{\theta=1}\\
 t&Q_{1/2}(t)&\|R(t)\|^2&Q_1(t)&\|R(t)\|^2\\\hline
20&0.00740737&0.00740741&0.000370342&0.000370370\\
40&0.00370369&0.00370370&0.0000925856&0.0000925926\\
80&0.00185184&0.00185185&0.0000231464&0.0000231481\\
160&0.000925922&0.000925926&0.00000578660&0.00000578704
\end{array}
\]
The relative gaps are at most $7.6\times10^{-5}$. The packet locations depend on $t$, as required by optimization in the source space. The script is provided with the supplementary files. These calculations concern the inherited-translation component. Theorem~\ref{thm:sharp-abstract-transfer} transfers the optimal exponents to the coupled system under its multiplier hypotheses.

\FloatBarrier

\section{Conclusion}
\label{sec:conclusion}

For systems with Hilbert-valued infinite memory, the decay from centered weighted histories is governed by the exact norm of the weighted right translation. Under the stated energy and multiplier hypotheses, and for densities of at most polynomial decay satisfying a logarithmic-derivative bound, the energy is bounded by the initial weighted norm times $\kappa(t)+\rho_w(t)$. When $\rho_w$ is integrable and the kernel has a finite first moment and a tail bound, or when $\kappa\le C_\kappa\rho_w$, the squared semigroup norm from $\mathscr X_w$ to $\mathscr X$ is bounded above and below by constant multiples of $\rho_w(t)$. For the kernel $pb(1+s)^{-p-1}$, with $p>1$, and the weight $(1+s)^{2\theta}$, with $\theta>0$, the exponent is $\min\{2\theta,p+1\}$ and saturates at $\theta=(p+1)/2$. For logarithmic kernels and weights, the logarithm of the weight governs the rate below this threshold and the logarithm of the kernel governs it above. Under the hypotheses of the transfer theorem, together with (K3) and a finite first moment, the tail weight $1+\kappa/m$ gives the exact rate $(1+t)^{-1}$. A prescribed-source theorem gives the sharp exponent $p-r$ for histories with two-sided algebraic growth, and it sharpens a borderline estimate of Guesmia for a polynomial subclass.

The weighted right-translation norm approach retains the inherited history as a decreasing source in the memory--dissipation splitting and closes the estimate by a linear dyadic recursion. It covers gradient-memory waves without instantaneous damping, endpoint-damped strings, damped plates and oscillators. The extension of the two-sided law to kernels with flat zones or jumps remains open.

\section*{Statements and declarations}

\noindent\textbf{Funding.} This research received no external funding.

\medskip
\noindent\textbf{Competing interests.} The author declares that there are no competing interests.

\end{document}